\documentclass{article}
\usepackage[usenames,dvipsnames]{color}

\usepackage{amsmath,amssymb,bm} 
\usepackage{arxiv}

\usepackage[utf8]{inputenc} % allow utf-8 input
\usepackage[T1]{fontenc}    % use 8-bit T1 fonts
\usepackage{hyperref}       % hyperlinks
\usepackage{url}            % simple URL typesetting
\usepackage{booktabs}       % professional-quality tables
\usepackage{amsfonts}       % blackboard math symbols
\usepackage{amsmath}        % mathematics enhanced
\usepackage{nicefrac}       % compact symbols for 1/2, etc.
\usepackage{microtype}      % microtypography
\usepackage{lipsum}
\usepackage{graphicx}
\usepackage{subcaption}
\usepackage{tikz,pgfplots,pgf} % THE package for drawing all sorts of things
\usetikzlibrary{matrix,shapes,arrows,positioning,decorations.pathreplacing}
\usetikzlibrary{plotmarks}

\usepackage{todonotes}

\newcommand{\intgt}{ {\cal T}(t) }

\graphicspath{ {./images/} }

\title{Cell dynamics in microfluidic devices under heterogeneous chemotaxis and growth conditions: a mathematical study}

\author{
 Jacobo Ayensa-Jiménez \\
  Aragon Institute of Engineering Research\\
  University of Zaragoza\\
  Mariano Esquillor, s.n. 50018, Zaragoza \\
  \texttt{jacoboaj@unizar.es} \\
   \And
 Mohamed H. Doweidar \\
  Mechanical Eng. Department, \\
  School of Engineering and Architecture (EINA) \\
  University of Zaragoza\\
  María de Luna s/n, Edificio Betancourt, 50018, Zaragoza \\
  \texttt{mohamed@unizar.es} \\
  \And
 Manuel Doblaré \\
  Aragon Institute of Engineering Research\\
  University of Zaragoza\\
  Mariano Esquillor, s.n. 50018, Zaragoza \\
  \texttt{mdoblare@unizar.es} \\
  \And
 Eamonn A. Gaffney \\
  Wolfson Centre for Mathematical Biology  \\ Mathematical Institute\\
  University of Oxford\\
   Woodstock Road, Oxford OX2 6GG, UK\\
  \texttt{gaffney@maths.ox.ac.uk} \\
}

\begin{document}
\maketitle

\begin{abstract}
As motivated by studies of cellular motility driven by  spatiotemporal chemotactic gradients in microdevices, we develop a framework for constructing approximate analytical solutions for the location, speed and cellular densities for cell chemotaxis waves in heterogeneous fields of chemoattractant from the underlying partial differential equation models. In particular, such chemotactic waves are not in general translationally invariant travelling waves, but possess a spatial variation that evolves in time, and may even may oscillate back and forth in time, according to the details of the chemotactic gradients. The analytical framework exploits the observation that unbiased cellular diffusive flux is typically small compared to chemotactic fluxes and is first developed and validated for a range of exemplar scenarios. The framework is subsequently applied to more complex models considering the full dynamics of the chemoattractant and how this may be driven and controlled within a microdevice by considering a range of boundary conditions. In particular, even though solutions cannot be constructed in all cases, a wide variety of scenarios can be considered analytically, firstly providing global insight into the important mechanisms and features of cell  motility in complex spatiotemporal fields of chemoattractant. Such analytical solutions also provide a means of rapid evaluation of model predictions, with the prospect of application in computationally demanding investigations relating theoretical models and experimental observation, such as Bayesian parameter estimation. 

\end{abstract}

% keywords can be removed

\section{Introduction}

Most biological processes integrate different cell populations, extracellular matrix (ECM) properties, chemotactic gradients and physical cues,  constituting  a  complex, dynamic and interactive microenvironment \cite{chen1997geometric,curtis1978control,schwarz2005physical,carter1967haptotaxis,lo2000cell}. Cells are also constantly adjusting   to accommodate  their surroundings,  particularly for homeostatically maintaining the intracellular and extracellular environment within  physiological constraints \cite{bray2000cell}. In response to the different chemical and physical external stimuli, cells can modify their shape, location, internal structure and genomic expression, as well as their capacity to proliferate, migrate, differentiate, produce ECM or other biochemical substances, changing, in turn, the surrounding medium as well as sending new signals to other cells \cite{mousavi20133d,kumar2013cell,mousavi2014computational}. This two-way interaction between cells and their environment is crucial in physiological processes such as embryogenesis, organ development, homeostasis, repair, and long-term evolution of tissues and organs among others, as well as in pathological processes such as atherosclerosis or cancer \cite{huang2005cell,hanahan2011hallmarks,nagelkerke2015mechanical,quail2013microenvironmental}. Furthermore, developing novel frameworks to investigate and elucidate  these mechanisms and interactions is key to developing novel therapeutic strategies aiming at promoting (blocking) desirable (undesirable) cellular behaviours \cite{mousavi2013cell}.

In particular, due to the underlying complexity, \emph{in vivo} research -- both in humans and in animals -- is impeded by the fact it is difficult to control and isolate effects. Thus, a simpler alternative is using \emph{in vitro} experiments. Nevertheless, the predictive power  currently available, whether \emph{in vivo} or \emph{in vitro}, is still poor, as demonstrated the continuous drop in the number of new drugs appearing annually, despite billion-dollar investments \cite{scannell2012diagnosing,nolan2007s}. Indeed, structural three-dimensionality is one of the most important characteristics of biological processes  \cite{edmondson2014three}, but \emph in vitro cells are mostly cultured in a traditional Petri dish (2D culture), where cell behaviour is dramatically different from   real tissues \cite{wang2014three}. Recently, microfluidics has arisen as a powerful tool to recreate the complex microenvironment that governs tumour dynamics \cite{sackmann2014present,bhatia2014microfluidic}. This technique allows the  reproduction of numerous important features that are lost in 2D cultures, as well as testing drugs in a much more reliable and efficient way \cite{bersini2014microfluidic, boussommier2016microfluidics, jeon2015human, zervantonakis2012three, wu2011induction}. 

In addition to such \emph{in vitro} models, mathematical \emph{in silico} models  are a powerful tool for dealing with many problems in physics, engineering, and biology. In particular, cell population evolution models based on transport partial differential equations (PDEs) have been widely used to study many biological processes, including cancer \cite{byrne2010dissecting,altrock2015mathematics}. For instance, tumour development is a key example of a highly dynamic and complex biological process that originates from external signals or stimuli modulated by the particular microenvironment. Furthermore, when a given treatment is applied (surgery, chemotherapy, radiotherapy, immunotherapy, hormones or a combination thereof), the tumour and its microenvironment undergo significant alterations. This leads tumour cells to proliferate and generate microenvironments that promote the death of surrounding cell types and the survival of tumour cells that are more adaptable and resistant. That is why, when modelling the enormous variety and complexity of a tumour and its microenvironment, the resulting differential equations are highly non-linear and strongly coupled \cite{kitano2002computational,bearer2009multiparameter,eils2013computational}. The numerical resolution of the equations, especially in the era of high performance computing, has been extensively utilised in the simulation of  ``what if'' scenarios and the study of effects and hypotheses in isolation, something that is often impossible to do with \emph{in vivo} and \emph{in vitro} models \cite{byrne2006modelling,katt2016vitro}. In turn, the construction and exploitation  of \emph{in silico} experiments is thus being increasingly used in the early stages of designing drugs and therapies against tumours.

A particular niche of interest is  Glioblastoma (GBM), the most common and aggressive primary brain tumour \cite{brat2012glioblastoma}, with  extensive studies dedicated to    mathematical modelling  its evolution \cite{hatzikirou2005mathematical},   reproducing aspects of GBM  histopathology \cite{bearer2009multiparameter} and incorporating the influence of  tumour microenvironment (TME) chemical and mechanical cues   \cite{kim2016role}. It has been demonstrated that GBM progression is extensively  controlled by the local oxygen concentrations and gradients \cite{hatzikirou2012go}, motivating many studies to incorporate  the role of oxygen gradients and hypoxia in tumour progression \cite{ayuso2017glioblastoma, martinez2012hypoxic,frieboes2006integrated}. Some models of GBM have  reproduced cell culture evolution under different experimental configurations \cite{ayensa2020mathematical}, using a \emph{go-or-grow} transition switch, governed  by  nonlinear activation functions for the chemotaxis and growth. Such studies therefore implicate the balance between cell migratory and proliferative activity, together with  their relation to the different TME stimuli, as playing a key role in GBM evolution.

Nevertheless, the complexity of the equations to be solved often require numerical simulations that are  impractical, due to the high computational cost, especially in the resolution of inverse problems such as parameter estimation, model selection, the design of experiments, sensitivity studies, model structural analysis and Uncertainty Quantification (UQ). Although many modern techniques as Reduced Order Models (ROM) and metamodels using Artificial Intelligence (AI) have been developed in recent years \cite{perez2021predicting}, the existence of analytical solutions, although approximate, provides key information to test and validate numerical algorithms, inform a mechanism based understanding across parameter space and to  allow initial predictions of Quantities of Interest (QoI),  such as travelling  fronts, equilibria, the ranges of variation in the solution across parameter space and parameter sensitivities, among others. Indeed, some works in the last years have focused on the use of these techniques for analysing GBM progression \cite{perez2011bright,gerlee2016travelling,stepien2018traveling}.

The interaction of cells with a chemoattractant leads to a type of Keller - Segel (K-S) model \cite{keller1971traveling}, which generally have a rich structure as reviewed by G. Arumugam and J. Tyagi \cite{arumugam2021keller}. One of the main interests concerning the K-S model is the existence and characteristics of travelling waves (see for instance \cite{xue2011travelling}). In this work, we explore the dynamics of cell populations under gradients of a chemotactic agent for one-dimensional problems. We move beyond travelling waves to investigate evolving wave solutions in the heterogeneous environments that are often found in microdevices and physiology. This general class of problems allows the treatment of a wide variety of situations related to the evolution of tumours, while the  analysis of the associated PDEs enables the quantification of histopathological characteristics, such as the spread of pseudopalisades and the response of the population to oscillatory stimuli. This knowledge can be used for the design of experiments, to speed up the characterisation processes of cell populations and to validate or rule out possible models. 

In particular, we are interested in modelling cell motility dynamics in microfluidic devices, which are experimental platforms that have been demonstrated to accurately recreate biomimetic physiological conditions \cite{shin2012microfluidic}, with application  in bioengineering and biomedical research \cite{sackmann2014present}. In many situations, the chemotactic agent concentration may be assumed as known, either because it can be directly measured, or because its concentration can be computed by solving a diffusion problem that is, to good approximation, decoupled from the cell population field. 

To proceed, we first describe the structure of the mathematical problem associated with the response of cell populations to chemotactic gradients, together with  the general assumptions and hypotheses about the underlying mechanisms. We derive pertinent features about the solution field, for instance that it possesses a migratory structure with a transition zone wavefront. We are also able to estimate the wave front evolution and the shape of the solution profile. In particular, we compute an analytical solution for  specific exemplar cases associated with specific relevant experimental situations. These include a constant spatial gradient of chemoattractant, together with temporal oscillations associated with a fluctuating source, quadratic profiles of chemoattractant, offering additional nonlinear features, and an exponential profile of chemoattractants corresponding to the diffusion from a localised source. We further apply the general results to the analysis of a cell culture microfluidic experiment, representing a slightly simplified version for an  {\it in vitro} model of GBM progression, as developed in \cite{ayensa2020mathematical}, showing how the methods presented here can generate analytical results for the simulation of microdevice representations for migratory tumour cell dynamics. 

\section{Methods}

\subsection{The model}

\begin{figure}[!htbp]
\centering
\begin{subfigure}{.49\textwidth}
  \centering
  \includegraphics[clip=true,trim=0pt 20pt 0pt 0pt,width=0.8\textwidth]{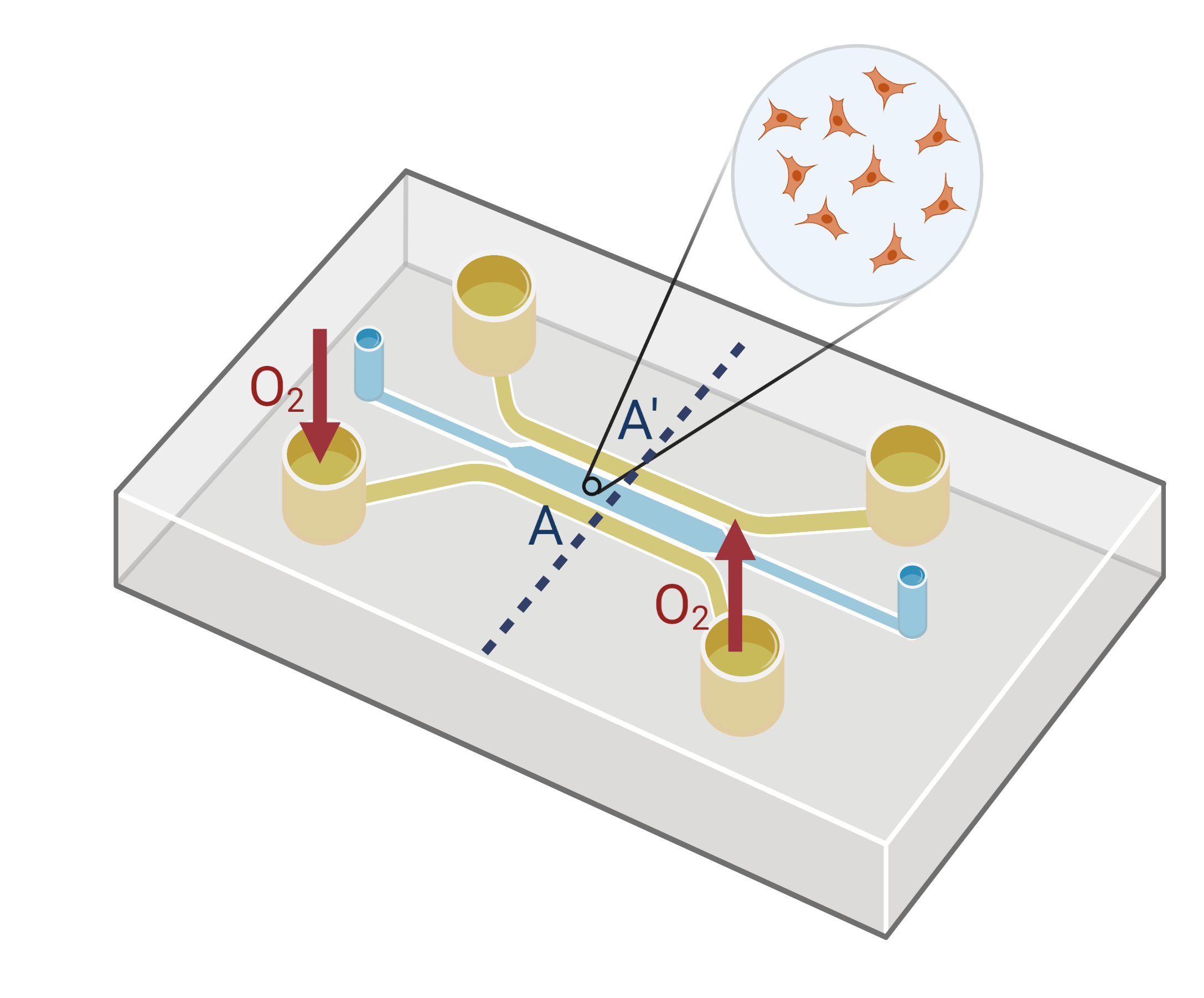}
  \caption{Scheme of the experimental configuration.}
  \label{fig:scheme_a}
\end{subfigure}
\begin{subfigure}{.49\textwidth}
  \centering
  \includegraphics[clip=true,trim=30pt 60pt 30pt 20pt,width=\textwidth]{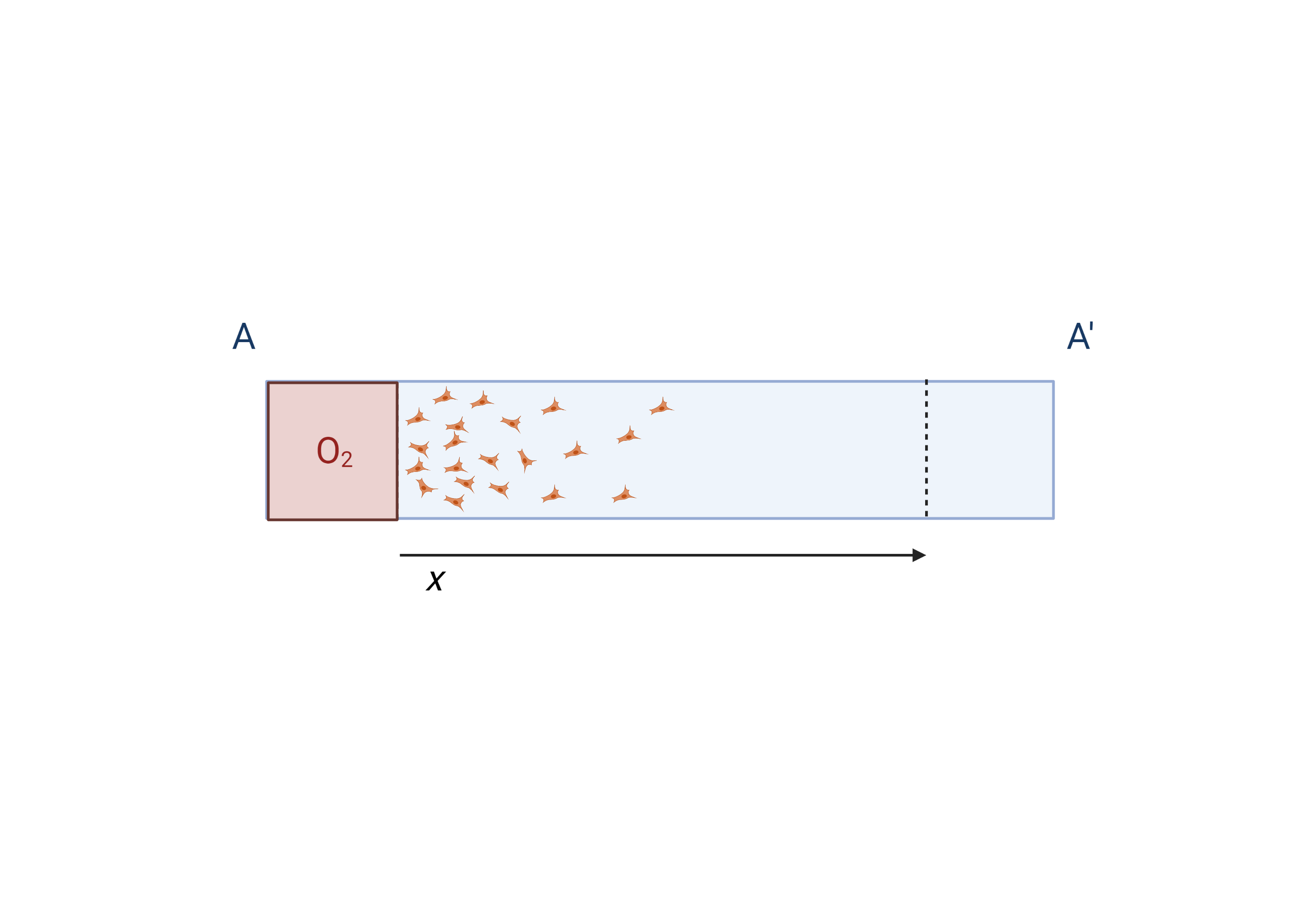}
  \caption{1D approximation of the cell culture.}
  \label{fig:scheme_b}
\end{subfigure}
\caption[Typical experimental configuration for modelling cell cultures]{\textbf{Typical experimental configuration for modelling cell cultures.} Due to the much larger length of the lateral channels relative to the width of the chamber, the domain geometry of the model is assumed one-dimensional, with length given by the length of the chamber, $L$. The cell concentrations are associated with a continuum field $u = u(x,t)$, with $t$ denoting time and $x$ the spatial coordinate along the chamber, as indicated. At the  lateral edges of the channel, that is $x=0,L$, zero flux boundary conditions are imposed, corresponding to the inability of cells to pass through these boundaries. Image created with BioRender.com.}
\label{fig:scheme}
\end{figure}

We study a broad class of problems that are related to the dynamics  of a cell culture in microfluidic devices under the influence of a chemotattractant, when the concentration of the agent can be computed or measured. A schematic view of this situation is represented in Fig. \ref{fig:scheme}.

The non-dimensional equation for the cell population concentration, $u$, represents cellular  diffusion and chemotaxis in response to a heterogeneous field of chemoattractant that  generates an advective flux  $\alpha(t,x)u$, and  also impacts cellular proliferation via $\beta(t,x)$ so as to  generate the governing equation
\begin{equation} \label{eq:eq_governing}
u_t + \left(\alpha(t,x)u\right)_x = D u_{xx} + \beta(t,x)u(1-u), 
\end{equation}
with $D>0$ the non-dimensional cellular diffusion coefficient. This governing equation is also supplemented by zero-flux boundary conditions, given by 
\begin{subequations} \label{eq:eq_boundary}
\begin{align}
    \left.Du_x-\alpha(t,x)u\right|_{x=0} &= 0, \\
    \left.Du_x-\alpha(t,x)u\right|_{x=L} &= 0,
\end{align}
\end{subequations}
and initial conditions
\begin{equation} \label{eq:eq_initial}
u(t=0,x)=u_0(x).
\end{equation}

\subsection{Computation of the general solution for small diffusion}

\subsubsection{Outer solution}

\begin{figure}[h]
\centering
\begin{tikzpicture}
\draw [black, ->] (0,0) -- (0,5);
\node [left] at (0,5) {$t$};
\draw [black, ->] (-0.3,0) -- (10,0);
\node [below] at (10,0) {$x$};
\node [left] at (-0.3,0) {$t_0$};
\draw [black, dashed] (0,2) -- (10,2);
\node [left] at (-0.3,2) {$t^*$};
\draw [black, dashed] (4,0) -- (4,4);
\draw (0,0) .. controls (2,3) and (3,0) .. (4,2);
\draw (4,2) .. controls (5,4) and (6,5) .. (8,4);

\draw [red, dashed] (0,1) .. controls (2,4) and (3,1) .. (4,3);
\draw [red, dashed] (4,3) .. controls (5,5) and (6,6) .. (8,5);

\draw [red, dashed] (0,2) .. controls (2,5) and (3,2) .. (4,4);
\draw [red, dashed] (4,4) .. controls (5,6) and (6,7) .. (8,6);

\draw [blue, dashed] (1,0) .. controls (3,2) and (3.5,0.5) .. (4,1.5);
\draw [blue, dashed] (4,1.5) .. controls (5,3) and (6,3) .. (8,3);

\draw [blue, dashed] (2,0) .. controls (3,1) and (3.5,0) .. (4,0.5);
\draw [blue, dashed] (4,0.5) .. controls (5,3) and (6,3) .. (8,2);

%\draw[blue,thick] (4,0) -- (8,4);
%\draw[pink, fill = pink] (4,0) circle [radius = 0.1];
%\draw[purple, fill = purple] (6,2) circle [radius = 0.1];
%\draw[green, fill = green] (5,0) circle [radius = 0.1];
%\draw[purple,decorate,decoration={brace,raise = 6pt,amplitude=5pt}] (6.5,2) -- (6.5,0) node[pos=0.5,right=10pt,purple]{$\Delta t_0$};
%\draw[purple,decorate,decoration={brace,raise = 6pt,amplitude=5pt}] (5.8,1) -- (5.8,0) node[pos=0.5,right=10pt,purple]{$\Delta t$};
%\draw[black,decorate,decoration={brace,raise = 6pt,amplitude=5pt}] (6,-0.5) -- (4,-0.5) node[pos=0.5,below=10pt,black]{$\Delta x$};
%\node [below] at (4,0) {$x_{i-1}$};
%\node [below] at (6,0) {$x_{i}$};
%\node [below, green] at (5,-0.1) {$L$};
\end{tikzpicture}
\caption[Projection of the characteristic curves]{\textbf{Projection of the characteristic curves.} The two families of characteristic curves are shown in blue and red.}
\label{fig:cplane}
\end{figure}
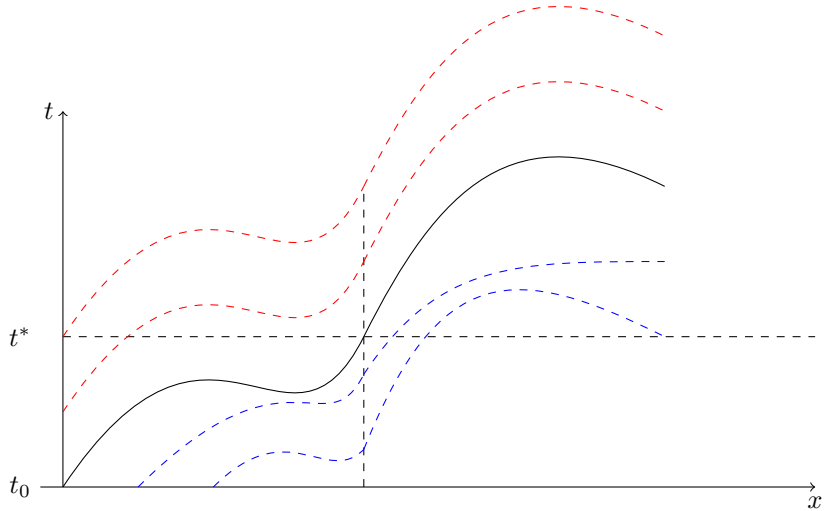

The main hypothesis, which is usually true for cellular motility due to weak cell-based random motility, is that the non-dimensional diffusion coefficient satisfies $D \ll 1$, as we will verify below in the case of GBM cells in microdevices. Hence, away from boundary layers,  diffusion may be neglected  compared to the influence of  growth and chemoattractant driven migration. Then,  Eq.~(\ref{eq:eq_governing}) may be approximated by:
\begin{equation} \label{eq:eq_governing_simp}
u_t + (\alpha(t,x)u)_x = \beta(t,x)u(1-u).
\end{equation}
This is a first-order hyperbolic PDE, amenable to  the method of characteristics. If we know the initial data  $u(0,x) = u_0(x)$, we can parameterise the initial data via  $s$ with the relation  $(t,x,u) = (0,s,u_0(s))$. Also, there is another family of characteristic curves, emerging from the $(t,x)$ points where $x=0$ and $t>0$, with the imposition of the $x=0$ boundary condition of no flux,  Eq.~(\ref{eq:eq_boundary}). Assuming that $\alpha(0,t) \neq 0$, and that the boundary is away from the transition region of the cellular wavefront, so that to excellent approximation $u_x= 0$ since spatial gradients are small,  we conclude from the boundary contition, Eq. (\ref{eq:eq_boundary}), that $u(0,t)=0$ to the same level of approximation. Therefore, this boundary condition can be parameterised via $s$ with the relation $(t,x,u) = (s,0,0)$. Hence, at $t=0$ and $x=0$ there is an emerging singular characteristic that splits the domain in two regions. The geometric interpretation of the method of characteristics is shown in Fig. \ref{fig:cplane}, which shows the projection of the characteristic curves onto the plane $(t,x)$.

We have, therefore, two families of characteristic curves. For the first one:
\begin{subequations} \label{eq:cF1_general}
\begin{align}
\frac{\mathrm{d}t}{\mathrm{d}\tau} &= 1, \quad t(0) = 0, \label{eq:cF1a_general} \\
\frac{\mathrm{d}x}{\mathrm{d}\tau} &= \alpha(t,x), \quad x(0) = s, \label{eq:cF1b_general}\\
\frac{\mathrm{d}u}{\mathrm{d}\tau} &= u\left(\beta(t,x)(1-u)-\alpha_x(t,x)\right), \quad u(0) = u_0(s), \label{eq:cF1c_general}
\end{align}
\end{subequations}
and for the second:
\begin{subequations} \label{eq:cF2_general}
\begin{align} 
\frac{\mathrm{d}t}{\mathrm{d}\tau} &= 1, \quad t(0) = s, \label{eq:cF2a_general} \\
\frac{\mathrm{d}x}{\mathrm{d}\tau} &= \alpha(t,x), \quad x(0) = 0, \label{eq:cF2b_general}\\
\frac{\mathrm{d}u}{\mathrm{d}\tau} &= u\left(\beta(t,x)(1-u)-\alpha_x(t,x)\right), \quad u(0) = 0. \label{eq:cF2c_general}
\end{align}
\end{subequations}

Noting the uniqueness of the solution to Eqs.~(\ref{eq:cF2_general}) courtesy of Picard's theorem,  we have by inspection that these equations only generate  the trivial solution $u(x,t)=0$.

However, the solution to the first family, 
given by Eqs.~(\ref{eq:cF1_general}), is typically more complex, though one always has \begin{equation} \label{eq:t_sol}
    t = \tau.
\end{equation}
Progress can be readily made when 
\begin{itemize} 
\item  $\alpha(t,x)$ is linear in $x$, such that $\alpha(t,x)=a(t)x+b(t)$ 
\item  $\alpha(t,x)$ is separable, that is $\alpha(t,x)=f(x)g(t)$. 
\end{itemize} 
In particular when  $\alpha(t,x)$ is linear in $x$ we have 
\begin{eqnarray} \label{sollin} x = s \mathrm{e}^{\int_0^t a(\eta)\mathrm{d} \eta} + 
\int_0^t b(\eta)\mathrm{e}^{\int_{\eta}^t a(\xi)\mathrm{d} \xi}   \, \mathrm{d} \eta =: F(t;s), 
\end{eqnarray}
which defines $F(t;s)$ for  $\alpha(t,x)=a(t)x+b(t)$. We generalise this definition so that $x=F(t;s)$ on the characteristic curve given by the value of $s$, and whenever this relation can be uniquely inverted for $s$, we write $s=G(t;x)$.

In contrast when $\alpha(t,x)=f(x)g(t)$ Eqs.~(\ref{eq:cF1a_general}) and (\ref{eq:cF1b_general}) may be integrated to obtain 
\begin{equation} \label{eq:x_sol}
    \int_s^x\frac{\mathrm{d}z}{f(z)} = \int_0^\tau g(\eta)\,\mathrm{d}\eta .
\end{equation}
In turn, Eq.~(\ref{eq:x_sol}) generates the relation  $x=F(t;s)$ on the characteristic curve. 
Further motivation for examples of linear and separable chemotactic response functions for $\alpha(t,x)$ are given in Section  \ref{sec:app} below. 

In both cases, or even in the more general case where $F(t;s)$ cannot readily be determined analytically for all relevant $t,s$, the location of the transition from $u=0$, and thus the location of the transition region for the cellular  wavefront, $x^* = x^*(t)$, is given by the characteristic with $s=0$ -- hence   $x^*(t)=F(t;0)$. This is particularly  informative about the general behaviour of the solution, for instance in determining the wavespeed. With respect to Eq.~(\ref{eq:cF1c_general}), we proceed using the change of variable $r = 1/u$, we obtain the ODE in terms of the $\tau = t$ variable:
\begin{equation} \label{eq:edo}
    r'(\tau)+\left(\beta(\tau,x(\tau;s))-H(\tau;s)\right)r =\beta(\tau,x(\tau;s)),
\end{equation}
where for the linear case $H(\tau; s) = a(\tau)$ and $H(\tau;s) = f'(F(\tau;s))g(\tau)$ for the separable case.

Noting the integration is along a characteristic, and thus  $s$ is fixed,  this equation is of the form $r'+p(\tau)r=q(\tau)$ for $p(\tau) =\beta(\tau,x(\tau;s))-H(\tau;s)$ and $q(\tau) = \beta(\tau,x(\tau;s))$, with $s$ fixed, so a general expression is given by
\begin{equation} \label{eq:edo_sol}
    r(\tau,s) = \exp\left(-\int_0^\tau p(\eta,s) \, \mathrm{d}\eta \right)\left[\frac{1}{u_0(s)} + \int_0^\tau q(\eta,s) \exp\left(\int_0^\eta p(\xi,s) \, \mathrm{d}\xi \right)\, \mathrm{d}\eta \right],  
\end{equation}
where $u_0(s)$ is indeed the value of $u_0$ at the location of the characteristic when $\tau=t=0$. 

Recapping, suppose $x = F(t;s)$ may be inverted to give $~s=G(t;x)$. Then, noting  $$ x(0,s)=F(0;s)=s,$$ by the parameterisation of the initial data, we have
$$ u(x,t) = u(x = F(t;s),t) = 1/r(t;s = G(t;x)),$$ and, in particular $$ u_0(s) = u_0(F(0,G(t;x))) = u_0(G(t;x)).$$ 
Combining these expressions with Eq.~(\ref{eq:edo_sol}), we obtain a general expression for $r$ and therefore for $u = u(x,t)$:
\begin{equation} \label{eq:u_sol_prev}
    u(x,t) = \left.\left( \frac{u_0(s)\exp\left(\int_0^t p(\eta,s) \, \mathrm{d}\eta \right)}{1 +  u_0(s)\int_0^t q(\eta,s) \exp\left(\int_0^\eta p(\xi,s) \, \mathrm{d}\xi \right)\, \mathrm{d}\eta}\right)\right|_{s = G(t;x)}, 
\end{equation}
where $s=G(t;x) $ is fixed on each characteristic curve and 
\begin{align}
    p(\eta,s) &= \beta(\eta,F(\eta;s))-H(\eta;s), \nonumber \\
    q(\eta,s) &= \beta(\eta,F(\eta;s)).  \nonumber 
\end{align}   
Eq. (\ref{eq:u_sol_prev}) may be also written in terms of $x$ and $t$ directly, obtaining
\begin{equation} \label{eq:u_sol}
    u(x,t) = \frac{u_0(G(t;x))\exp\left(\int_0^t p(\eta,G(t;x)) \, \mathrm{d}\eta \right)}{1 +  u_0(G(t;x))\int_0^t q(\eta,G(t;x)) \exp\left(\int_0^\eta p(\xi,G(t;x)) \, \mathrm{d}\xi \right)\, \mathrm{d}\eta}. 
\end{equation}

%\subsubsection{A special case for \texorpdfstring{$\bm{u(x,t).}$}.}
\paragraph{A special separable case for \texorpdfstring{$\bm{u(x,t).}$}.}

Often below, when the chemoattractant flux term is independent of time, so that $\alpha(t,x)=f(x)$, we will have $u_0\ll 1$, $\beta(t,x)\equiv 1$, $g(t)=1$. In these circumstances we have the simplification
$$s=G(t;x)=F(-t;x),$$
by the symmetry $(x,s,\tau)\rightarrow(s,x,-\tau)$ in Eq.~(\ref{eq:x_sol}) with $g(t)=1$ 
and  also that 
$$p(\eta,s) = 1 - f'(F(\eta;s)), $$ which gives 
$$u(x,t) =  u_0(G(t;x))\mathrm{e}^t \exp\left(-\int_0^t f'(F(\eta,s)) \, \mathrm{d}\eta \right) +\mathcal{O}(u_0^2). $$
Recalling that the integration is along a characteristic, so that $s$ is fixed,  we change the integration variable via $X= F(\eta;s)$, noting from Eq.  (\ref{eq:t_sol}) and from Eq.~(\ref{eq:cF1b_general}), with $s$ fixed and $\alpha(t,x) = f(x)$, that
$$ \frac 1 {f(X(\eta))}= \frac{\mathrm{d}\eta}{\mathrm{d}X}.$$ Hence, on further noting $s=F(0;s),~x=F(t;s)$ and $s = G(t,x)$, we have, to within $\mathcal{O}(u_0^2)$ corrections, that 
\begin{equation} 
u(x,t) = u_0(G(t;x)) \mathrm{e}^t \exp\left(-\int_{s}^{x} \frac{f'(X)}{f(X)} \, \mathrm{d}X\right)  = u_0(G(t;x)) \mathrm{e}^t \frac{f(G(t;x))}{f(x)}   = u_0(G(t;x)) \mathrm{e}^t \frac{f(F(-t;x))}{f(x)}. \label{u0approx}
\end{equation}

In the examples plotted below we also have $u_0(x)=u_0$ is constant, in which case $  u_0(G(t;x))$ collapses to the constant $u_0$. For clarity please note  that $u_0$ is an abbreviation for  $ u_0(G(t;x))$ in general, though this  is  constant and denoted simply by $u_0$ in all examples, separable or otherwise, plotted below.

\subsubsection{Inner solution}

To explore the transition layer moving with the wavefront, we introduce a scaling of coordinates such that diffusion and advection provide a leading order dominant balance in the transition layer, with $x=x^*(t)+\delta X$ for $X$ the inner variable and $\delta \ll 1 $. With the change of variables 
$$ (t,x)\rightarrow(\tau, X), ~~~ ~~~  t=\delta \tau,  ~~~~~~  x=x^*(t)+\delta X, ~~~~~~  U( \tau,X)=u(t,x),$$
 one has the inner solution equations 
\begin{eqnarray*} \frac{1}{\delta}U_\tau &=&\frac 1 \delta \left( x_t^*-\alpha(t,x^*(t)+\delta X ) \right)U_X + \frac{D}{\delta^2}u_{XX} +\beta U(1-U) - \alpha_x(t,x^*(t)+\delta X) \\ &=&   \frac 1 \delta \left( \alpha(t,x^*(t) ) -\alpha(t,x^*(t)+\delta X ) \right)U_X + \frac{D}{\delta^2}U_{XX} +\mathcal{O}(1),  
\end{eqnarray*}
where the second line uses $x_t^*=\alpha(t,x^*(t) )$. We take  $\delta =D \ll 1 $ to  bring the advective and diffusive terms into a nominal dominant balance.  

%We cannot assume the transition layer is time-independent, firstly as it evolves in time  and also because neglecting the temporal dynamics ultimately leads to a breakdown of the perturbation method. Hence  we rescale time via $\tau=\delta \tau^*$ to bring the time derivative into the nominal dominant balance. 
%With a further, weak, assumption that $\alpha $ uniformly possesses an order one derivative with respect to $x$, one obtains  

Let us now assume that $\alpha$ uniformly possesses an order one derivative with respect to $x$, that is, we assume that $\alpha(t,x^*(t)) - \alpha(t,x^*(t) + \delta X) \simeq -\alpha_x(t,x^*(t))\delta X$, therefore
$$ U_{\tau} = u_{XX} +\mathcal{O}(\delta),$$ 
for $X\sim \mathcal{O}(1)$. The solution of this equation, when $U(X,\tau=0) = H(X)$, with $H$ the Heaviside step function is
$$U(X,\tau) =  \frac{1}{2}\left(1+\mathrm{erf}\left(\frac{X}{2\sqrt{\tau}}\right)\right).$$

%At leading order, thus dropping the $\mathcal{O}(\delta)$ terms, this has an error function solution when the initial conditions for the diffusion equation are a Heaviside.  This is self-consistent though we  do we need the numerics to confirm that the fast initial transients, too fast to be seen in the plots below, do generate an approximate Heaviside after a very short amount of time -- this is clear from the plots in Figs 3--6 at $t=0.04$. 

\subsubsection{Composite solution}

We can combine the inner solution, rewritten in terms of $(x,t)$ with the outer solutions via a composite approximation to generate an approximation for the full numerical solution across the domain (away from any prospective boundary layer at $x=1$). Thus, with $u_0(x)=u_0$, constant, and $u_{\mathrm{an}}(x,t)$ the analytical characteristic solution of Eq \eqref{eq:u_sol} we have 
\begin{equation}
u(x,t) \sim
\left\{
 \begin{array}{ccc}
   u_{\mathrm{an}}(x^*(t),t)U\left(\frac{x-x^*(t)}{D},\frac{t}{\tau}\right) + &0, & ~~ x< x^*(t), \\
   u_{\mathrm{an}}(x^*(t),t)U\left(\frac{x-x^*(t)}{D},\frac{t}{\tau}\right) + &\left(u_\mathrm{an}(x,t) -  u_\mathrm{an}(x^*(t),t)\right), & ~~ x\geq x^*(t),
 \end{array} \right.
\end{equation}
so, finally:
\begin{equation}
u(x,t) \sim
\left\{
 \begin{array}{cc}
   \frac{1}{2}u_{\mathrm{an}}(x^*(t),t)\left(\mathrm{erf}\left(\frac{x-x^*(t)}{2\sqrt{tD}}\right)+1\right), & ~~ x< x^*(t), \\
   \frac{1}{2}u_{\mathrm{an}}(x^*(t),t)\left(\mathrm{erf}\left(\frac{x-x^*(t)}{2\sqrt{tD}}\right)-1\right) + u_\mathrm{an}(x,t), & ~~ x\geq x^*(t),
 \end{array} \right.
\end{equation}

For the special case discussed in the previous section,  that is, when $\alpha(t,x) = f(x)$, $\beta(t,x) \equiv 1$ and $u_0 \ll 1$, $u_\mathrm{an}(x,t) = u_0(G(t;x))\mathrm{e}^t\frac{f(G(t;x))}{f(x)}$ and in particular $u_\mathrm{an}(x^*(t),t) = u_0(0)\mathrm{e}^t\frac{f(0)}{f(x)}$ so that for this special case 
\begin{equation}
u(x,t) \sim
\left\{
 \begin{array}{cc}
   \frac{1}{2}u_0(0)\mathrm{e}^t\frac{f(0)}{f(x)}\left(\mathrm{erf}\left(\frac{x-x^*(t)}{2\sqrt{tD}}\right)+1\right), & ~~ x< x^*(t), \\
   \frac{1}{2}u_0(0)\mathrm{e}^t\frac{f(0)}{f(x)}\left(\mathrm{erf}\left(\frac{x-x^*(t)}{2\sqrt{tD}}\right)-1\right) +  u_0(G(t;x))\mathrm{e}^t\frac{f(G(t;x))}{f(x)} , & ~~ x\geq x^*(t),
 \end{array} \right.
\end{equation}
with continuity assured from the constraint $G(t;x^*(t))=0$, which holds as, by construction,  both $s=0$ and $x=x^*(t)$ on the separating characteristic and $s=G(t;x)$.

\subsection{Some particular cases of interest}\label{spcoi}

The solution given by Eq.~(\ref{eq:u_sol}) is general given that   $D \ll 1$, though even in the separable case using Eq.~(\ref{eq:x_sol}) to determine the functions $F$ and $G$ can generate complicated solutions that are not readily expressed in terms of standard functions.  Further, even when $\alpha(t,x)$ is separable or linear, allowing extensive analytical progress, there is still considerable freedom in the form of $\alpha(t,x)$ and hence we firstly analyse models where $\alpha(t,x)$ is linear, quadratic, and exponential in space, before proceeding to consider an example of cellular behaviour in a  microdevice. 

\begin{figure}[!h]
    \centering
    \includegraphics[clip=true,trim=0pt 0pt 0pt 0pt,width=0.7\textwidth]{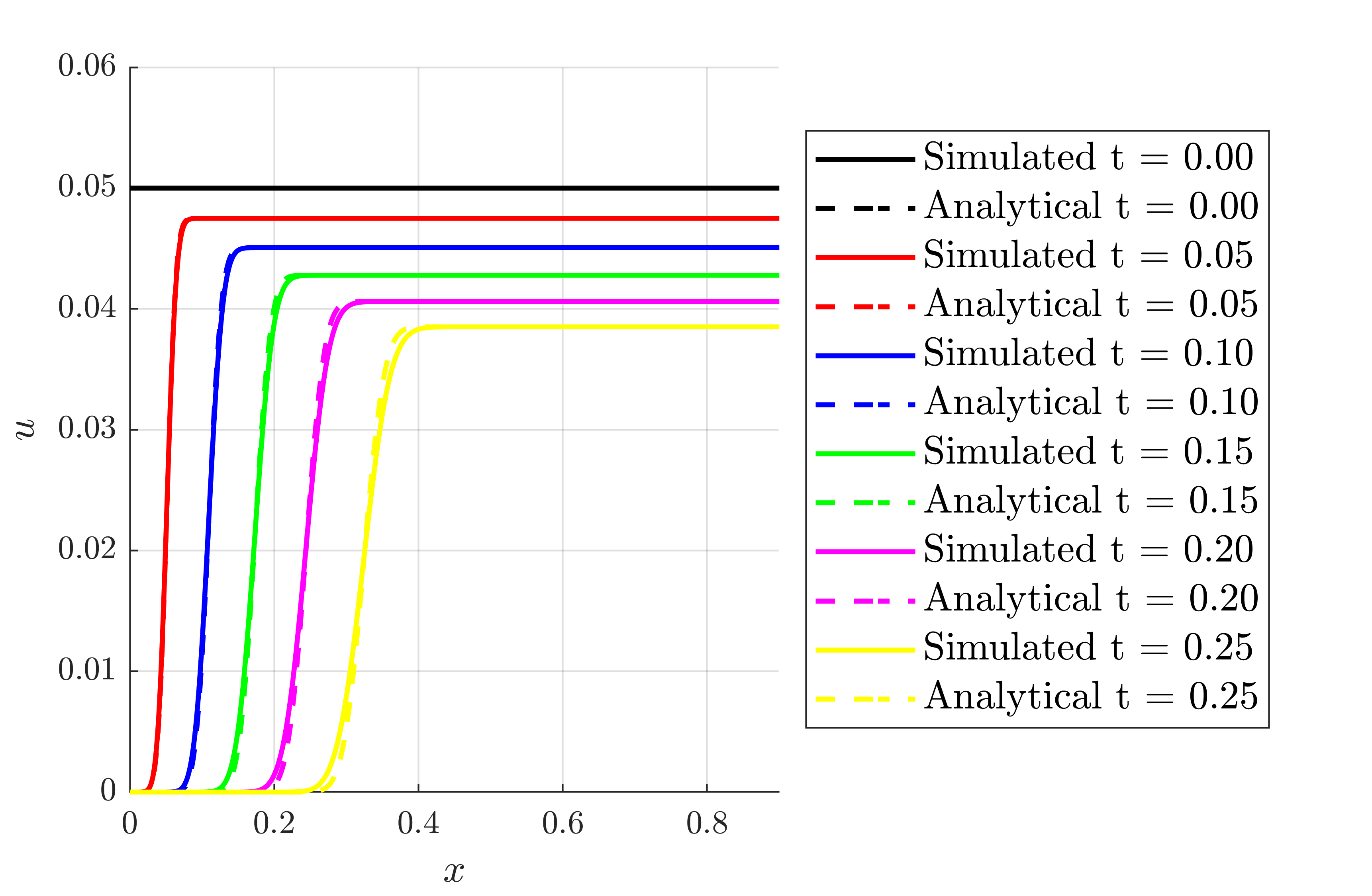}
    \caption[Comparison of numerical and analytical solutions for $\alpha$ linear.]{\textbf{Comparison of numerical and analytical solutions for the linear case.} The analytical and simulated profiles are compared at different times, considering $a = 2$ and $b = 1$, $\beta(t,x) = g(t) = 1$ and $u_0(x) = u_0 = 0.05$ whereas the full numerical simulation was obtained using $D=1\times 10^{-3}$.  For this plot, and for analogous plots below, there are localised boundary layer effects near the right hand edge of the domain, $x=1$, that are not captured by the presented analytical solution.}\label{fig:lin_comp}
\end{figure}

\subsubsection{Linear chemotaxis} \label{sec:linear_chemo}

We first consider a function of the form 
\begin{equation}
    \alpha(t,x) = (ax+b)g(t).
\end{equation}

Then, we have $f(x)=ax+b$ and Eq.~(\ref{eq:x_sol}) yields 
\begin{equation}
\frac 1 a \ln\left(\frac{ax+b}{as+b}\right) = \int_0^\tau g(\eta)\,\mathrm{d}\eta =: {\cal T}(\tau). 
\end{equation}
  Thus defines the function ${\cal T}$, and with ${\cal T}(t) = \int_0^t g(\eta)\,\mathrm{d}\eta $ we have 
\begin{subequations} \label{eq:FG_lin}
\begin{align}
    F(t;s) &= \frac{1}{a}\left[(as+b)\mathrm{e}^{a\intgt }-b\right], \nonumber \\
    G(t;x) &= \frac{1}{a}\left[(ax+b)\mathrm{e}^{-a\intgt }-b\right]. \nonumber
\end{align}
\end{subequations}
    
Hence, the transition is located at
\begin{equation} \label{eq:xstar_lin}
    x^*(t) = \frac{b}{a}\left(\mathrm{e}^{a\intgt  }- 1\right).
\end{equation}
    
Furthermore, for $\beta(t,x) = 1$ and $g(t)=1$, so that ${\cal T}(t)=t$,  the integral expression for $u(x,t)$ in Eq. \eqref{eq:u_sol} is readily determined to reveal  
$$ u(x,t) = \frac{ u_0 (a-1) \mathrm{e}^{-(a-1)t}}{a-1+u_0 \left(1- \mathrm{e}^{-(a-1)t}  \right)} =  u_0  \mathrm{e}^{-(a-1)t} +\mathcal{O}(u_0^2),$$ while $u=0$ before the transition.
Fig. \ref{fig:lin_comp} shows a comparison between the numerical results with $D=1\times 10^{-3}$ and the approximate analytical solution for $\alpha(x)=2x+1$ and $\beta(t,x) = 1 = g(t)$.  For this plot, and for analogous plots below, note the accurate prediction of the evolving front, $x^*(t)$ and the general agreement between the numerical and analytical solutions for $u(x,t)$.  Furthermore, one can expect boundary layer effects near the right hand edge of the domain, $x=1$, that are not captured by the presented solution, though these are not plotted in the current Figure.

A further, and particularly relevant case, is when $g(t)=\cos (\omega t)$, whence 
\begin{subequations} \label{eq:FG_osc}
\begin{align}
    F(t;s) &= \frac{(as+b)\exp\left(\frac{a}{\omega}\sin(\omega t)\right)-b}{a}, \\
    G(t;x) &= \frac{(ax+b)\exp\left(-\frac{a}{\omega}\sin(\omega t)\right)-b}{a}.
\end{align}
\end{subequations}
Note we still have  $G(t;x) = F(-t,x)$ since $g(t)=\cos(\omega t)$ is even so that $(x,s,\tau)\rightarrow(s,x,-\tau)$ remains  a symmetry of  Eq.~(\ref{eq:x_sol}). 
The transition is located at
\begin{equation} \label{eq:xstar_osc}
    x^*(t) = \frac{b}{a}\left[\exp\left(\frac{a}{\omega}\sin(\omega t)\right)-1\right].
\end{equation}
    
\begin{figure}[!t]
    \centering
    \includegraphics[clip=true,trim=0pt 0pt 0pt 0pt,width=0.7\textwidth]{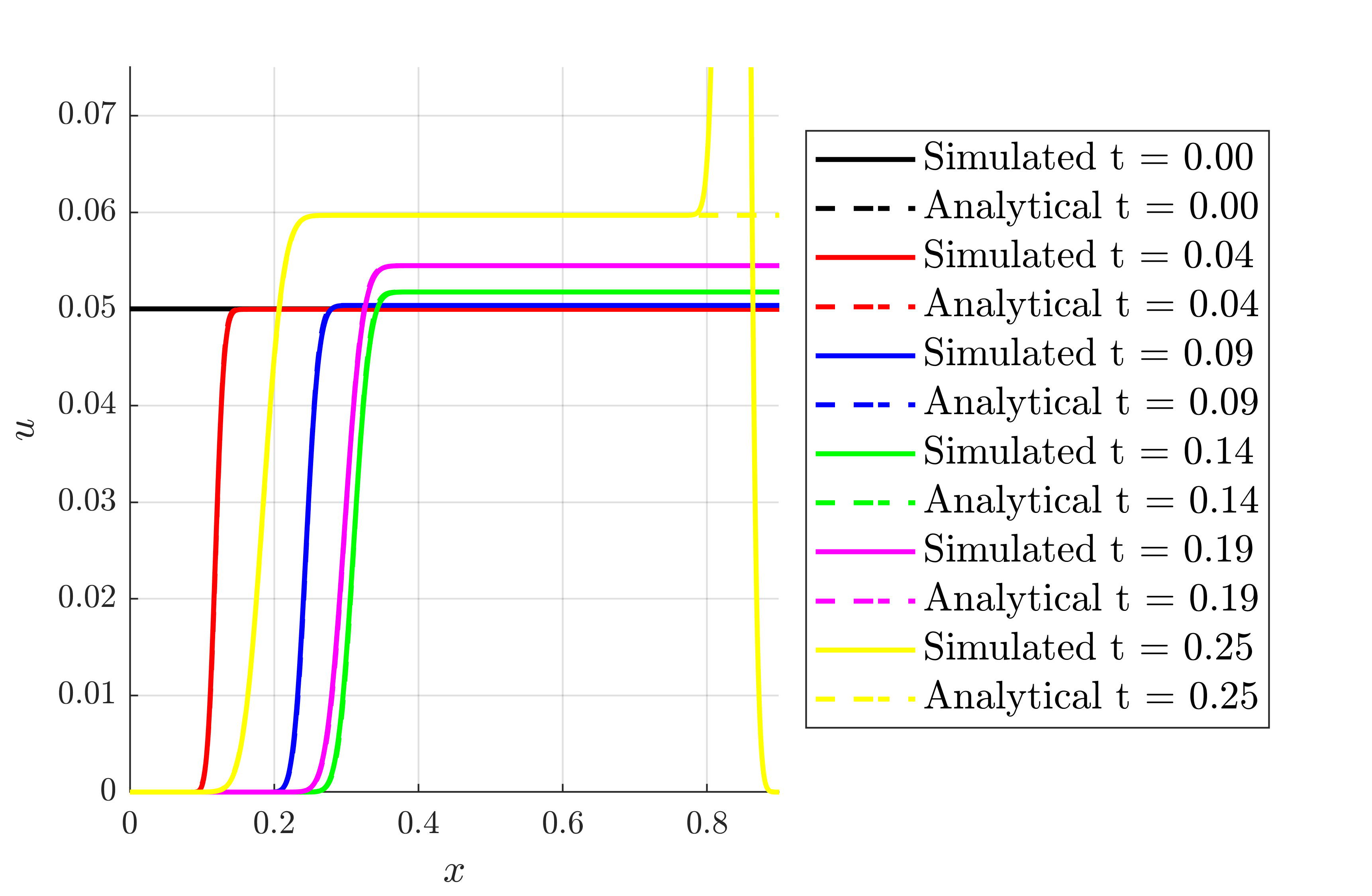}
    \caption{{\bf Comparison between analytical solutions using asymptotic theory and numerical solutions for the oscillatory gradient.}  The analytical and simulated profiles at different times are compared, considering $a=1$, $b=3$, $\omega = 10$, $\beta(t,x) = 1$ and $u_0(x) = u_0 = 0.05$ whereas the full numerical simulation was again obtained using $D=1\times 10^{-3}$. One can clearly observe that the wave of cells oscillates, with the cell density spatially constant on the right of the transition, except on approaching $x=1$, where there is a boundary layer that is not accommodated in the analysis. \label{figosc}} \label{fig:osc_comp}
\end{figure}

Fig. \ref{fig:osc_comp} shows the comparison between the numerical results with $D=1\times 10^{-3}$ and the analytical solution, with $\beta(t,x) = 1$. Note that as all the domain is displayed, the boundary layer in the numerical solution near $x=1$ can be readily observed.

The oscillating gradient case is very pertinent  as it corresponds to the case where the oxygenation feed between the two channels at the microfluidic device switches, so we shall explore it in more detail. Rather than working with the general solution given by Eq. (\ref{eq:u_sol}) it can be more expedient to consider the differential equation  for $r=1/u$ given by Eq.~(\ref{eq:edo}); noting $\beta=1$, $f(x)=(ax+b)$, $u_0(x)=u_0$, constant, for the cases considered here, this reduces to  
\begin{equation} \label{eq:edo_osc_prev}
    r'+\left(1-a\cos(\omega t)\right) r = 1, \quad r(0) = 1/u_0, 
\end{equation}
with no $x$-dependence. Hence, on the right of the transition region the solution is constant, $u(x) = u^*$.

%We have obtained the general solution:
% \begin{equation} \label{eq:u_sol_osc}
%     u(x,t) = \left.\frac{u_0(G(t;x))\exp\left(\int_0^t p(\eta,s) \, \mathrm{d}\eta \right)}{1 + u_0(G(t;x))\int_0^t q(\eta,s) \exp\left(\int_0^\eta p(\xi,s) \, \mathrm{d}\xi \right)\, \mathrm{d}\eta}\right|_{s = G(t;x)},
% \end{equation}
%where, for the oscillating gradient case:
%\begin{align} 
%    p(\eta,s) &= \beta(F(\eta,s),\eta)-a\cos(\omega \eta), \nonumber \\
%    q(\xi,s) &= \beta(F(\xi;s),\xi), \nonumber \\
%    s &= G(t,x), \nonumber
%\end{align}
%and $F$ and $G$ are given by the expressions Eqs. (\ref{eq:FG_osc}). For some specific regimes, it is possible to derive approximate solutions to Eq. (\ref{eq:edo}), that is, more manageable expressions than Eq. (\ref{eq:u_sol_osc}), which is expressed by a quadrature and therefore may be %somewhat complex, for example, in inverse problems.

\begin{figure}
    \begin{subfigure}{0.49\textwidth}
        \centering
        \includegraphics[clip=true,trim=0pt 0pt 0pt 0pt,width=0.7\textwidth]{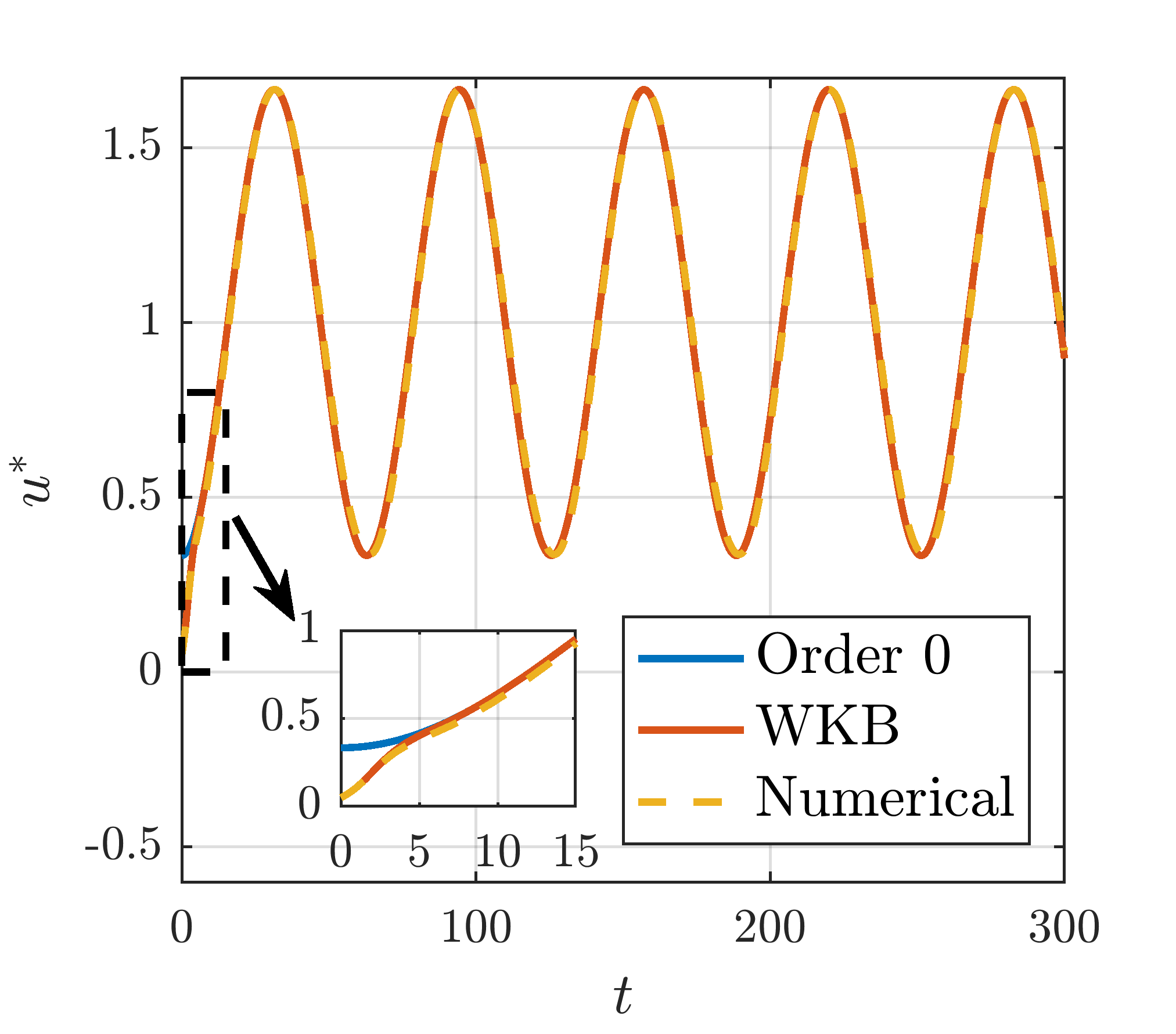}
        \caption{ $\omega \ll 1$ ($a =2$, $\beta = 3$ and $\omega = 0.1$).}\label{fig:omega_small}
    \end{subfigure} 
    \begin{subfigure}{0.49\textwidth}
        \centering
        \vspace{0.35cm}
        \includegraphics[clip=true,trim=0pt 0pt 0pt 0pt,width=0.7\textwidth]{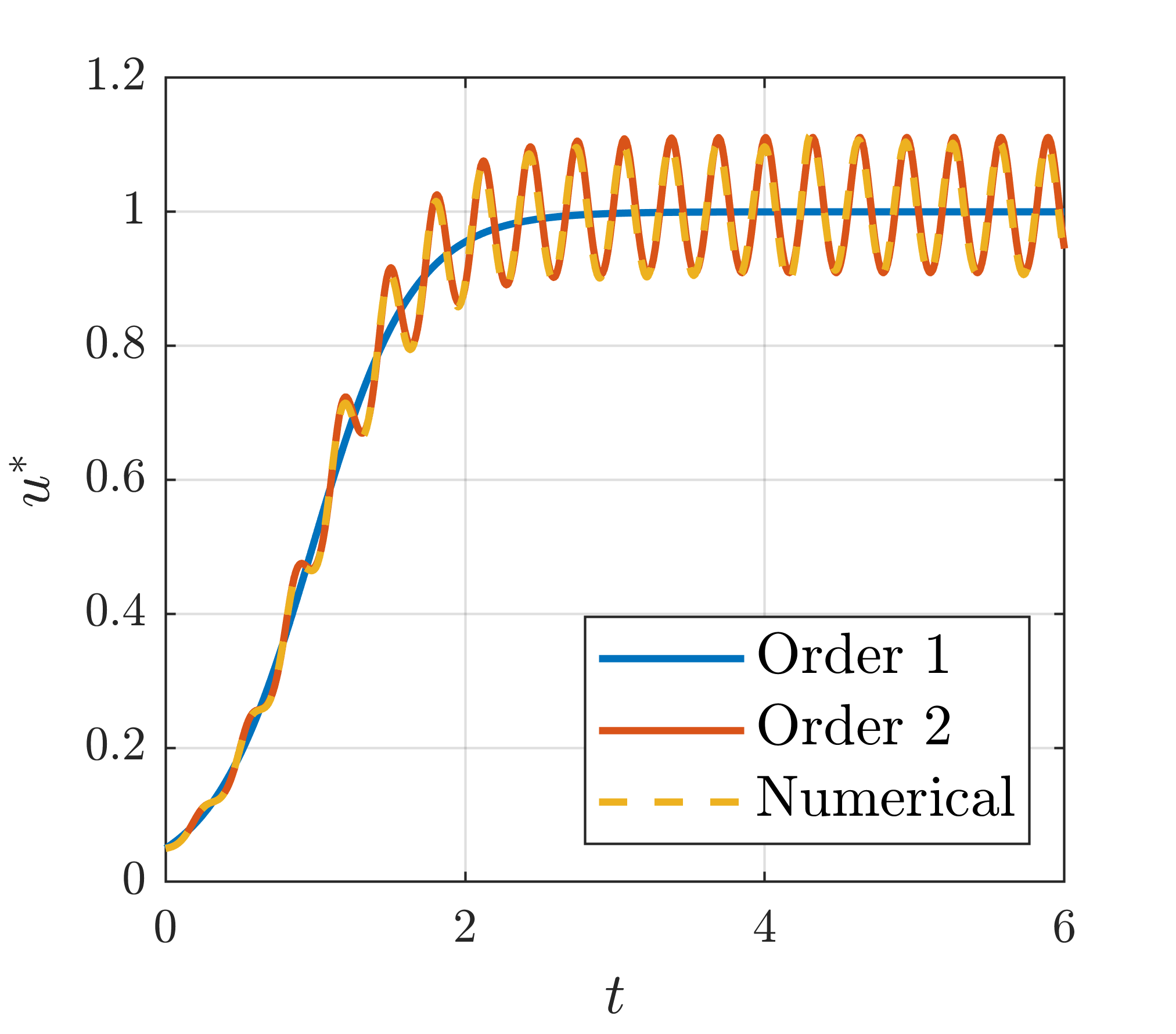}
        \caption{ $\omega \gg 1$ ($a =2$, $\beta = 3$ and $\omega = 20$).}\label{fig:omega_large}
    \end{subfigure} \\
    \begin{subfigure}{0.49\textwidth}
        \centering
        \includegraphics[clip=true,trim=0pt 0pt 0pt 0pt,width=0.7\textwidth]{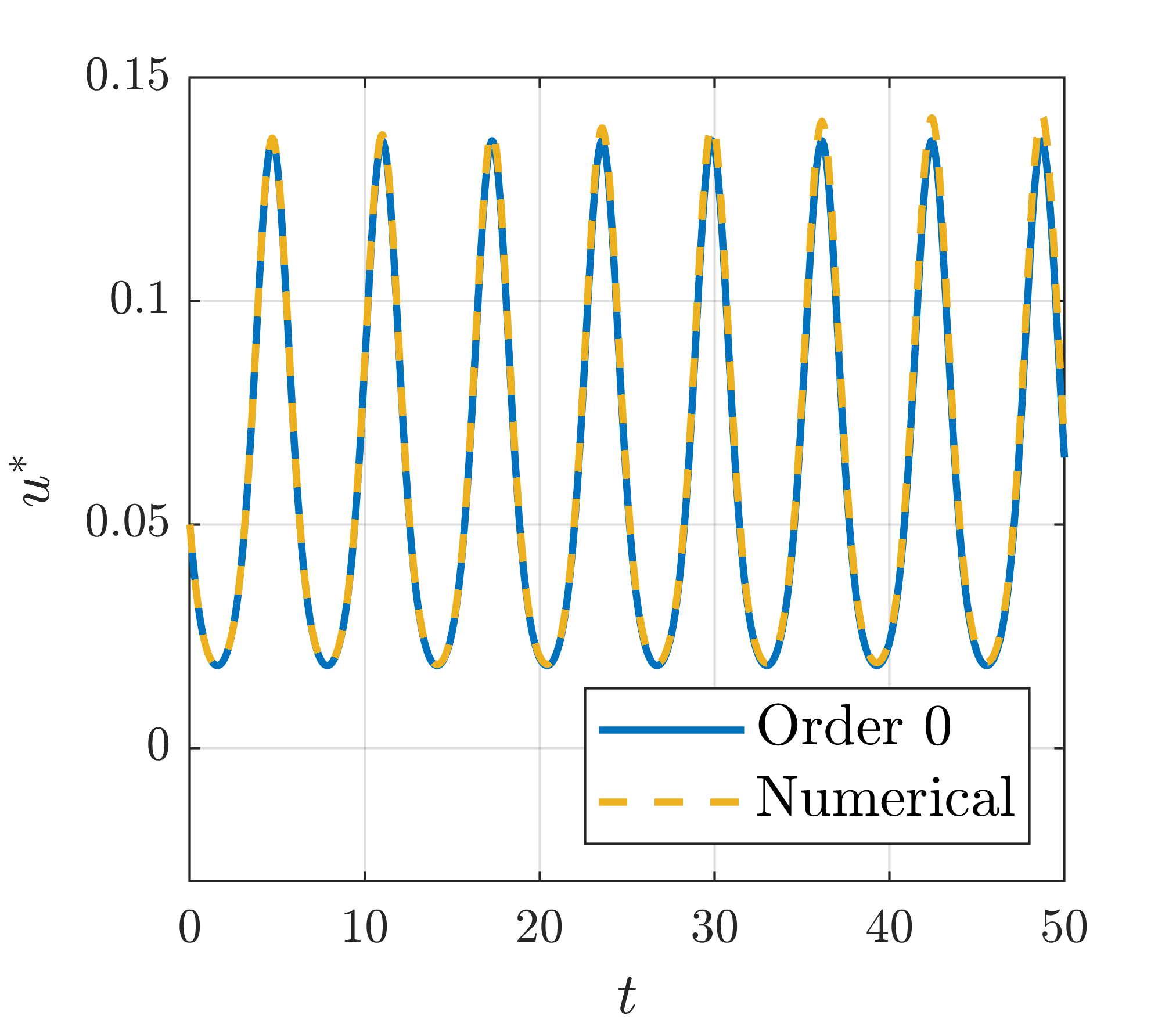}
        \caption{ $\beta \ll a$ ($a =1$, $\beta = 0.001$ and $\omega = 1$).}\label{fig:beta_small}
    \end{subfigure}
    \begin{subfigure}{0.49\textwidth}
        \centering
        \includegraphics[clip=true,trim=0pt 0pt 0pt 0pt,width=0.7\textwidth]{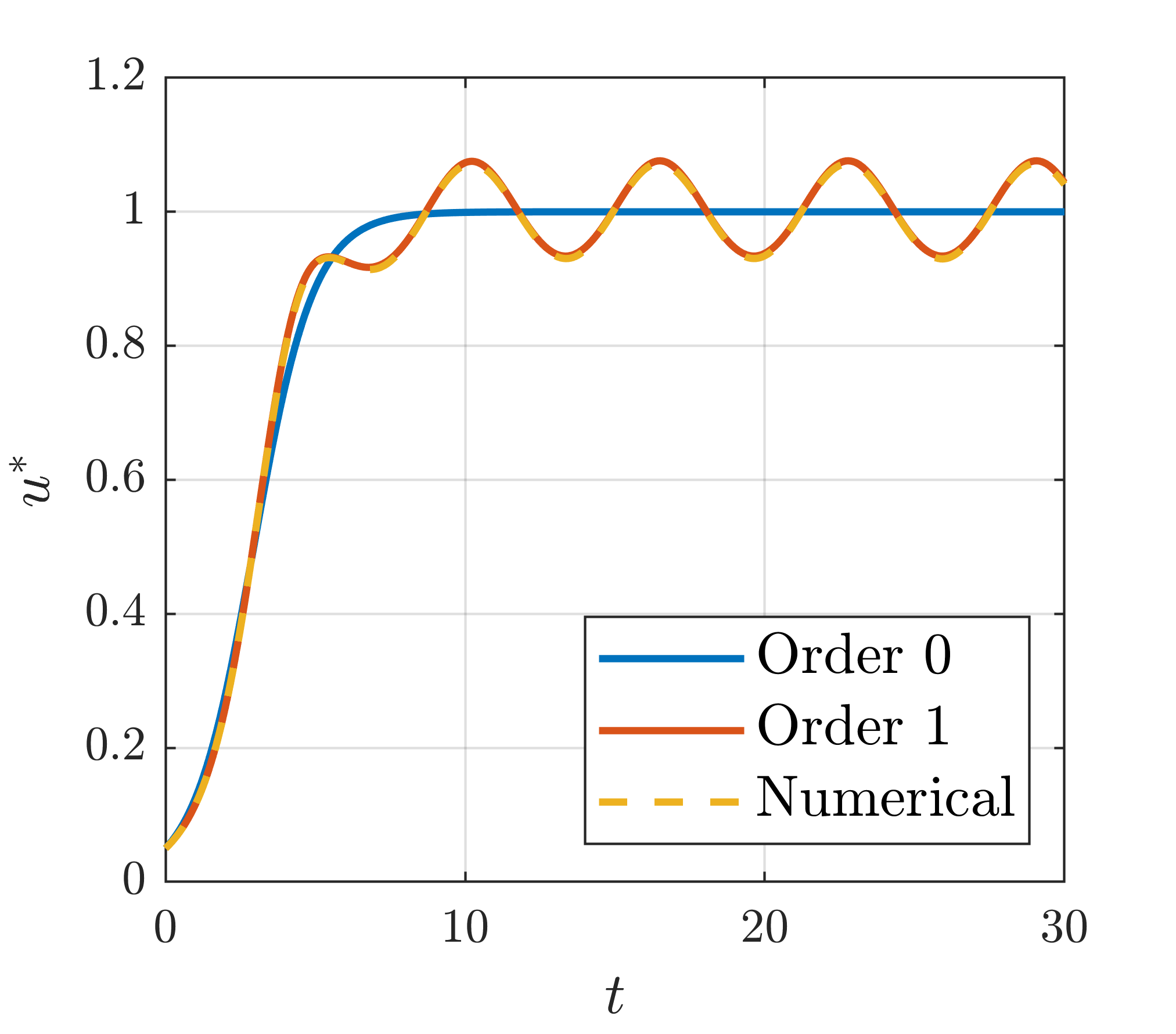}
        \caption{ $a \ll \beta$ ($a =0.01$, $\beta = 1$ and $\omega = 1$).}\label{fig:A_small}
    \end{subfigure}
    \caption{\textbf{Comparison between analytical solutions using asymptotic theory and numerical solutions for the oscillatory gradient}. The four exposed cases are analysed: slow variations of the gradiens ($\omega \ll 1$), fast variations of the gradients ($\omega \gg 1$), dominant chemotaxis ($\beta \ll a$) and dominant growth ($a \ll \beta$).}
    \label{fig:at_comp}
\end{figure}

In  Appendix \ref{sec:a1} we show that  Eq. \ref{eq:edo_osc_prev} may be solved using different asymptotic methods for four different regimes:
\begin{itemize}
    \item Slow variations of the gradients, $\omega \ll 1$. \vspace*{-1mm}
    \item Fast variations of the gradients, $\omega \gg 1$. \vspace*{-1mm}
    \item Dominant chemotaxis, $\beta \ll a$. \vspace*{-1mm}
    \item Dominant growth, $a \ll \beta$.\vspace*{-1mm}
\end{itemize}
For spatial locations to the right of the transition region, but away from any boundary layer at $x=1$, these approximate solutions are compared to  numerical solutions computed using standard Runge-Kutta solvers in Fig \ref{fig:at_comp}.

\clearpage

\subsubsection{Quadratic chemotaxis}

\begin{figure}[!p]
    \centering
    \begin{subfigure}[t]{\textwidth}
    	\centering
        \includegraphics[clip=true,trim=0pt 0pt 0pt 0pt,width=0.7\textwidth]{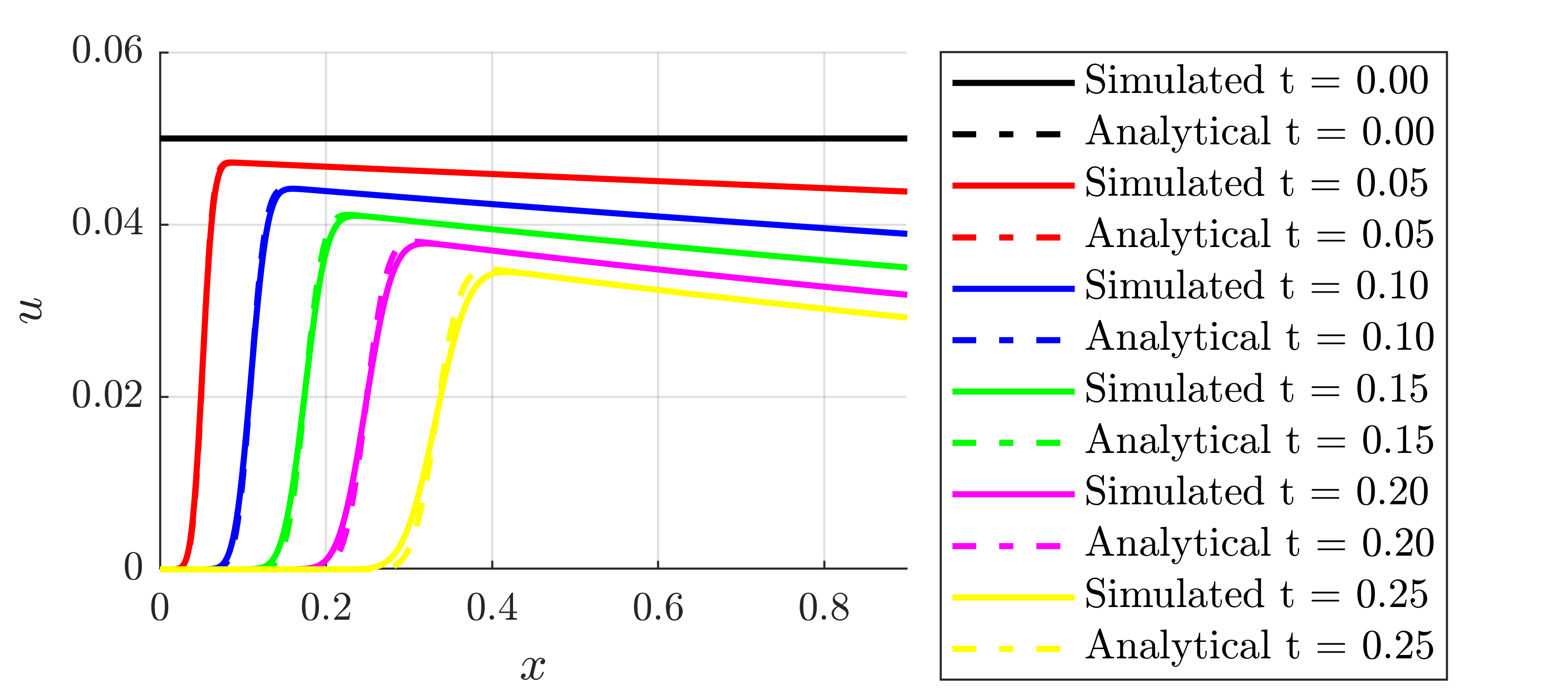}
        \caption{$a=1$, $b=2$, $c=1$ ($\Delta = 0$).}
        \label{fig:comp_poly_1}
    \end{subfigure} \\
    \begin{subfigure}[t]{\textwidth}
    	\centering
        \includegraphics[clip=true,trim=0pt 0pt 0pt 0pt,width=0.7\textwidth]{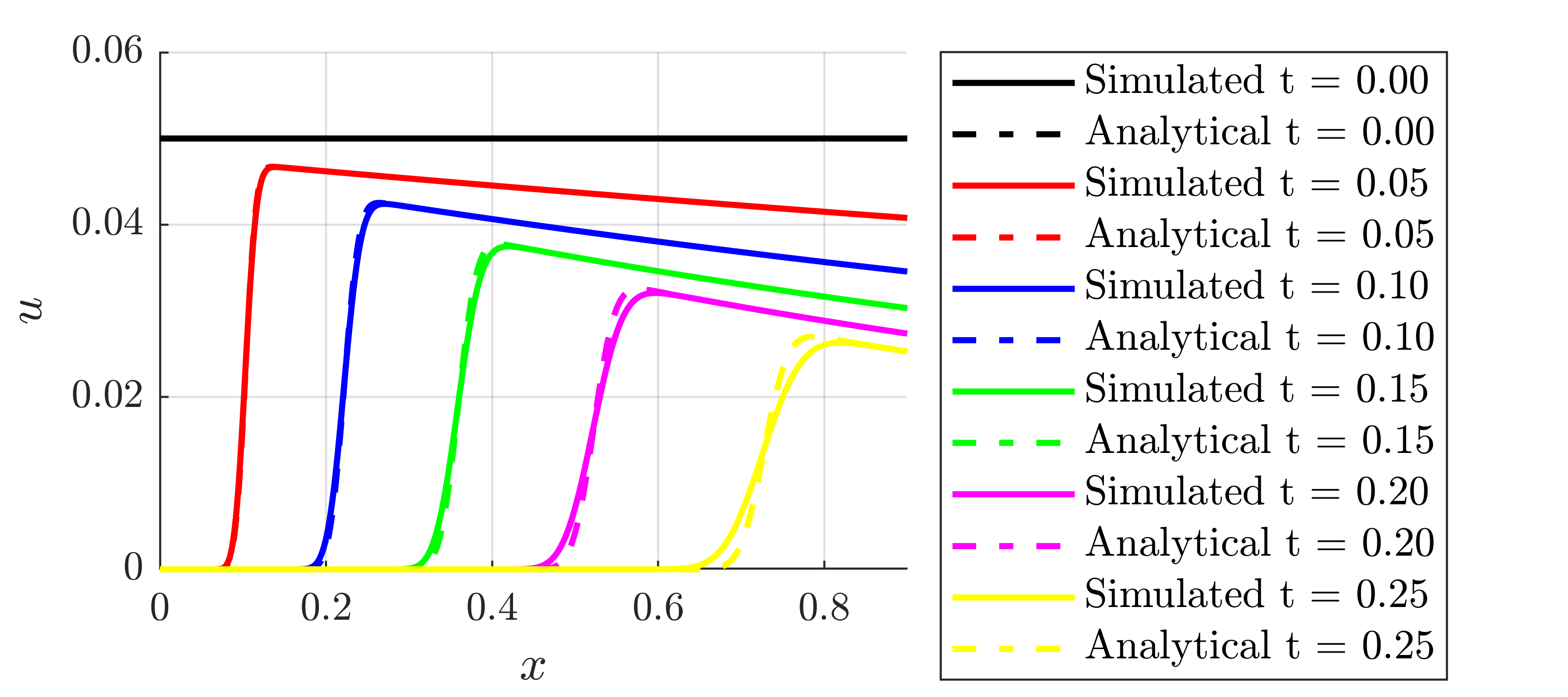}
        \caption{$a=2$, $b=2$, $c=2$ ($\Delta < 0$).}
        \label{fig:comp_poly_2}
    \end{subfigure} \\
    \begin{subfigure}[t]{\textwidth}
    	\centering
        \includegraphics[clip=true,trim=0pt 0pt 0pt 0pt,width=0.7\textwidth]{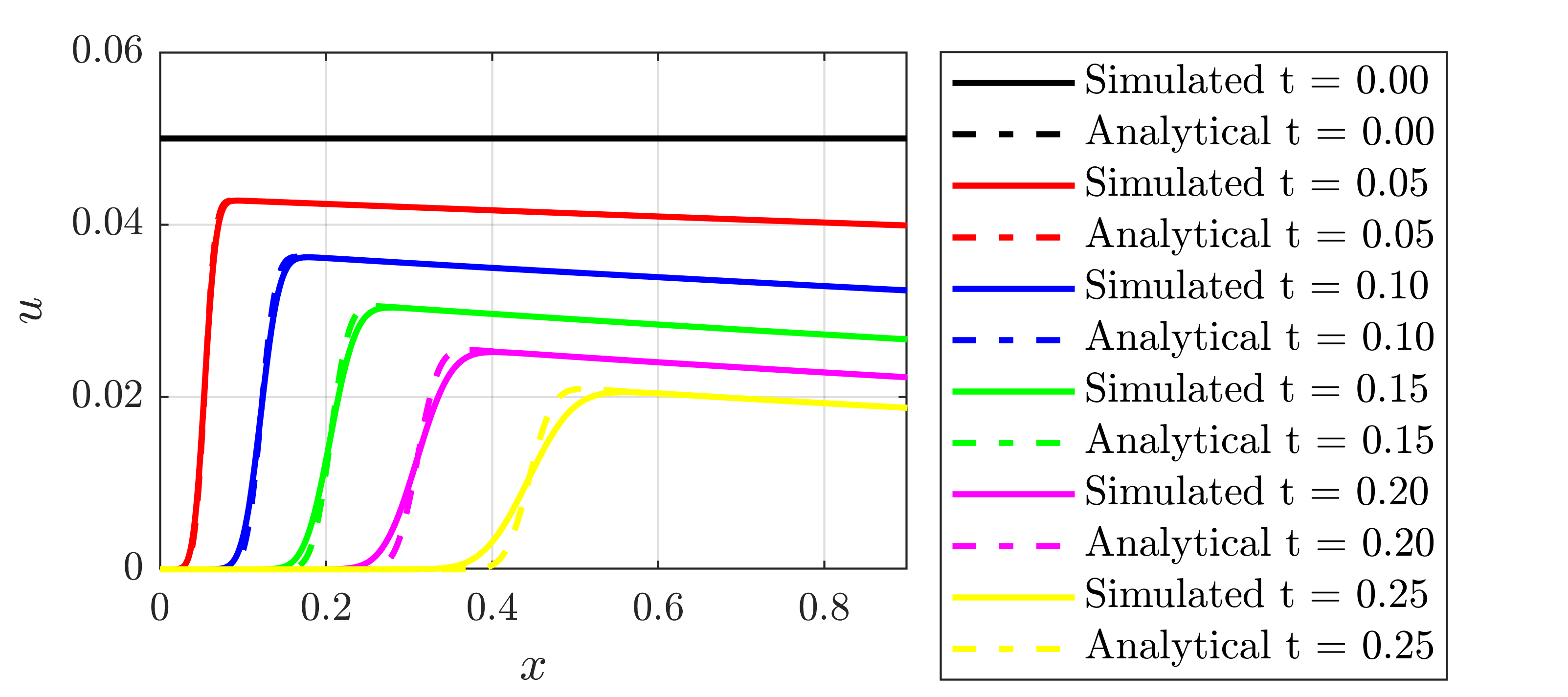}
        \caption{$a=1$, $b=4$, $c=1$ ($\Delta > 0$).}
        \label{fig:comp_poly_3}
    \end{subfigure}
    \caption[Comparison of numerical and analytical solutions for $\alpha$ quadratic.]{\textbf{Comparison of numerical and analytical solutions for $\alpha$ quadratic.} The analytical and simulated (for $D = 1 \times 10^{-3}$) profiles at different times are compared for three different $\alpha(t,x) = ax^2+bx+c$ expressions.}\label{fig:polynomial_comp}
\end{figure}

Now, we consider a function of the form 
\begin{equation}
    \alpha(t,x) = (ax^2+bx+c)g(t), 
\end{equation}
and hence $f(x)=ax^2+bx+c$.
\begin{itemize}
    \item If $\Delta = b^2-4ac = 0$, Eq.~(\ref{eq:x_sol}) gives 
    
    \begin{equation}
    \frac{1}{2as+b}-\frac{1}{2ax+b} = \frac{1}{2}\int_0^\tau g(\eta)\,\mathrm{d}\eta = \frac 1 2 {\cal T}(\tau),
    \end{equation}
    and in turn 
    \begin{subequations} \label{eq:FG_q1}
    \begin{align}
        F(t;s) &= \frac{4as+2abs\intgt +b^2\intgt }{4a-2ba\intgt -4a^2s\intgt }, \\
        G(t;x) &= \frac{4ax-2ab\intgt x-b^2\intgt }{4a^2 x\intgt +2ab\intgt +4a}. 
    \end{align}
    \end{subequations}
    
    The transition is located at 
    \begin{equation} \label{eq:xstar_q1}
        x^*(t) = \frac{b^2\intgt }{4a-2ba\intgt }  = \frac{c\intgt}{1-b\intgt/2}.
    \end{equation}
    Furthermore, in the case where $g(t)=\beta(t,x)=1$ we have 
    to an accuracy of $\mathcal{O}(u_0^2)$ that  
    $$ u(x,t) = u_0 \mathrm{e}^t \frac{f(s)}{f(x)} =  u_0 \mathrm{e}^t \left. \frac{as^2+bs+c}{ax^2+bx+c}\right|_{s=G(t;x)}= \frac{4u_0\mathrm{e}^t}{(2+bt +2atx)^2},$$
    via use of  Eq. (\ref{u0approx}), with $u_0\ll1$ where significant, but elementary, manipulation is required to deduce the final expression. 
    In particular for the parameters of Fig. \ref{fig:comp_poly_1} it is straightforward to show that at leading order the behaviour in $x$ for $t$ fixed is linear, as observed in this figure.

    \item If $\Delta = b^2-4ac < 0$, Eq.~(\ref{eq:x_sol}) yields 
    \begin{equation}
    \left(\arctan\left(\frac{2ax+b}{\sqrt{ -\Delta}}\right)-\arctan\left(\frac{2as+b}{\sqrt{ -\Delta}}\right)\right) = \frac{1}{2}\sqrt{ -\Delta}\int_0^\tau g(\eta)\,\mathrm{d}\eta = \frac{1}{2}\sqrt{ -\Delta}{\cal T}(\tau) ,
    \end{equation}
    so that 
    \begin{subequations} \label{eq:FG_q2}
    \begin{align}
        F(t;s) &= \frac{1}{2a}\left[\sqrt{-\Delta}\tan\left(\frac{1}{2}\sqrt{-\Delta}\intgt +\arctan\frac{2as+b}{\sqrt{-\Delta}}\right)-b\right], \\
        G(t;x) &= \frac{1}{2a}\left[\sqrt{-\Delta}\tan\left(-\frac{1}{2}\sqrt{-\Delta}\intgt +\arctan\frac{2ax+b}{\sqrt{-\Delta}}\right)-b\right], 
    \end{align}
    \end{subequations}
    with the transition is located at
    \begin{equation} \label{eq:xstar_q2}
        x^*(t) = \frac{1}{2a}\left[\sqrt{-\Delta}\tan\left(\frac{1}{2}\sqrt{-\Delta}\intgt +\arctan\frac{b}{\sqrt{-\Delta}}\right)-b\right].
    \end{equation}
   
      Furthermore, the analytical and full numeric solutions are plotted for  the parameters $a=b=c=2$ with $\beta(t,x)=g(t)=1$ in Fig. \ref{fig:comp_poly_2}. For this set of parameters,  where $t\leq 0.25$ as in Fig. \ref{fig:comp_poly_2}, we have  
   \begin{eqnarray} \tan \left(\frac{t\sqrt{-\Delta}}2\right) \approx \frac{t\sqrt{-\Delta}}2 \left(1+ O\left(\frac 1 3 \left( \frac{t\sqrt{-\Delta}}2  \right)^2\right)  \right) \approx  \frac{t\sqrt{-\Delta}}2,  ~~~~ \mbox{noting} ~~~~  \frac 1 3 \left( \frac{t\sqrt{-\Delta}}2  \right)^2 \lesssim 0.06 ~. \label{app1} 
   \end{eqnarray}
    Dropping the  corrections in higher powers of  $t\sqrt{-\Delta}$, valid for sufficiently small time including the times plotted in Fig. \ref{fig:polynomial_comp}, 
    we have 
    $$ x^*(t) =  \frac{ct}{1-bt/2} ,$$
    and then also dropping terms scaling with $\mathcal{O}(u_0^2)$,  one finds
      $$ u(x,t) = u_0 \mathrm{e}^t \frac{f(s)}{f(x)} =  u_0 \mathrm{e}^t \left. \frac{as^2+bs+c}{ax^2+bx+c}\right|_{s=G(t;x)}= \frac{u_0\mathrm{e}^t (4-t^2\Delta ) }{(2+bt +2atx)^2}, $$
      with a relative correction of $(1+\mathcal{O}([t\sqrt{-\Delta}/2] ^2/3)).$
The latter  again gives an approximate  linear dependence in $x$ for $t$ fixed given the parameters of  Fig. \ref{fig:comp_poly_2}, as observed. 
    
    \item If $\Delta = b^2-4ac > 0$, Eq.~(\ref{eq:x_sol}) gives 
    \begin{equation}
    \ln\left(\frac{2ax+b-\sqrt{\Delta}}{2ax+b+\sqrt{\Delta}}\right)-\ln\left(\frac{2as+b-\sqrt{\Delta}}{2as+b+\sqrt{\Delta}}\right) = \sqrt{\Delta}\int_0^\tau g(\eta)\,\mathrm{d}\eta = \sqrt{\Delta}\mathcal{T}(\tau),
    \end{equation}
    so that 
    \begin{subequations} \label{eq:FG_q3}
    \begin{align}  F(t;s) & = \frac{1}{2a}\left[\frac{\gamma^+(2as + \gamma^-)\exp(\frac{1}{2}\sqrt{\Delta}\intgt) - \gamma^-(2as + \gamma^+)\exp(-\frac{1}{2}\sqrt{\Delta}\intgt)}
        {(2as+\gamma^+)\exp(-\frac{1}{2}\sqrt{\Delta}\intgt) - (2as+\gamma^-)\exp(\frac{1}{2}\sqrt{\Delta}\intgt)}\right],  \\
     G(t;x) & = \frac{1}{2a}\left[\frac{\gamma^+(2ax + \gamma^-)\exp(-\frac{1}{2}\sqrt{\Delta}\intgt) - \gamma^-(2ax + \gamma^+)\exp(\frac{1}{2}\sqrt{\Delta}\intgt)}
        {(2ax+\gamma^+)\exp(\frac{1}{2}\sqrt{\Delta}\intgt) - (2ax+\gamma^-)\exp(-\frac{1}{2}\sqrt{\Delta}\intgt)}\right], 
    \end{align}
    \end{subequations}
    where we have defined $\gamma^+ = b+\sqrt{\Delta}$ and $\gamma^- = b-\sqrt{\Delta}$, with the transition location given by 
    \begin{equation} \label{eq:xstar_q3}
        x^*(t) = 2c\left[\frac{\exp(\frac{1}{2}\sqrt{\Delta}\intgt) - \exp(-\frac{1}{2}\sqrt{\Delta}\intgt)}{\gamma^+\exp( -\frac{1}{2}\sqrt{\Delta}\intgt) - \gamma^-\exp( \frac{1}{2}\sqrt{\Delta}\intgt})\right].
    \end{equation}
    
     For  Fig \ref{fig:comp_poly_3}, we again have $\beta(t,x)=1=g(t)$, so that ${\cal T}(t)=t$; we also have $t\sqrt{\Delta}/2 \ll 1$ throughout the simulation regime. Hence, on neglecting higher powers of $t\sqrt{-\Delta}/2 $, the transition location simplifies to 
    $$x^*(t) = \frac{ct}{1- bt/2 },$$
   which agrees with the above as $t\sqrt{-\Delta} / 2 \rightarrow 0.$ Furthermore, under these conditions with $u_0^2\ll 1$ and  neglecting higher orders in $u_0$, one finds  
     $$ u(x,t) = u_0 \mathrm{e}^t \frac{f(s)}{f(x)} =  u_0 \mathrm{e}^t \left. \frac{as^2+bs+c}{ax^2+bx+c}\right|_{s=G(t;x)}=  u_0 \mathrm{e}^t(1-(b+2ax)t) , $$
      with an relative error 
     of $(1+\mathcal{O}([t\sqrt{-\Delta}/2]^2)).$  For $t$ fixed, again we have $u(x,t)$ is approximately linearly decreasing in $x$, as observed in Fig  \ref{fig:comp_poly_3}.
\end{itemize}

More generally, Fig. \ref{fig:polynomial_comp} shows a comparison between the numerical results with $D=1\times 10^{-3}$, $\beta(t,x)=1=g(t)$ and the analytical solution for three different expressions of $\alpha(x)$ and $\beta(t,x) = 1$, with full consistency with the above solutions and approximations.

\begin{figure}[!t]
    \centering
    \includegraphics[clip=true,trim=0pt 0pt 0pt 0pt,width=0.7\textwidth]{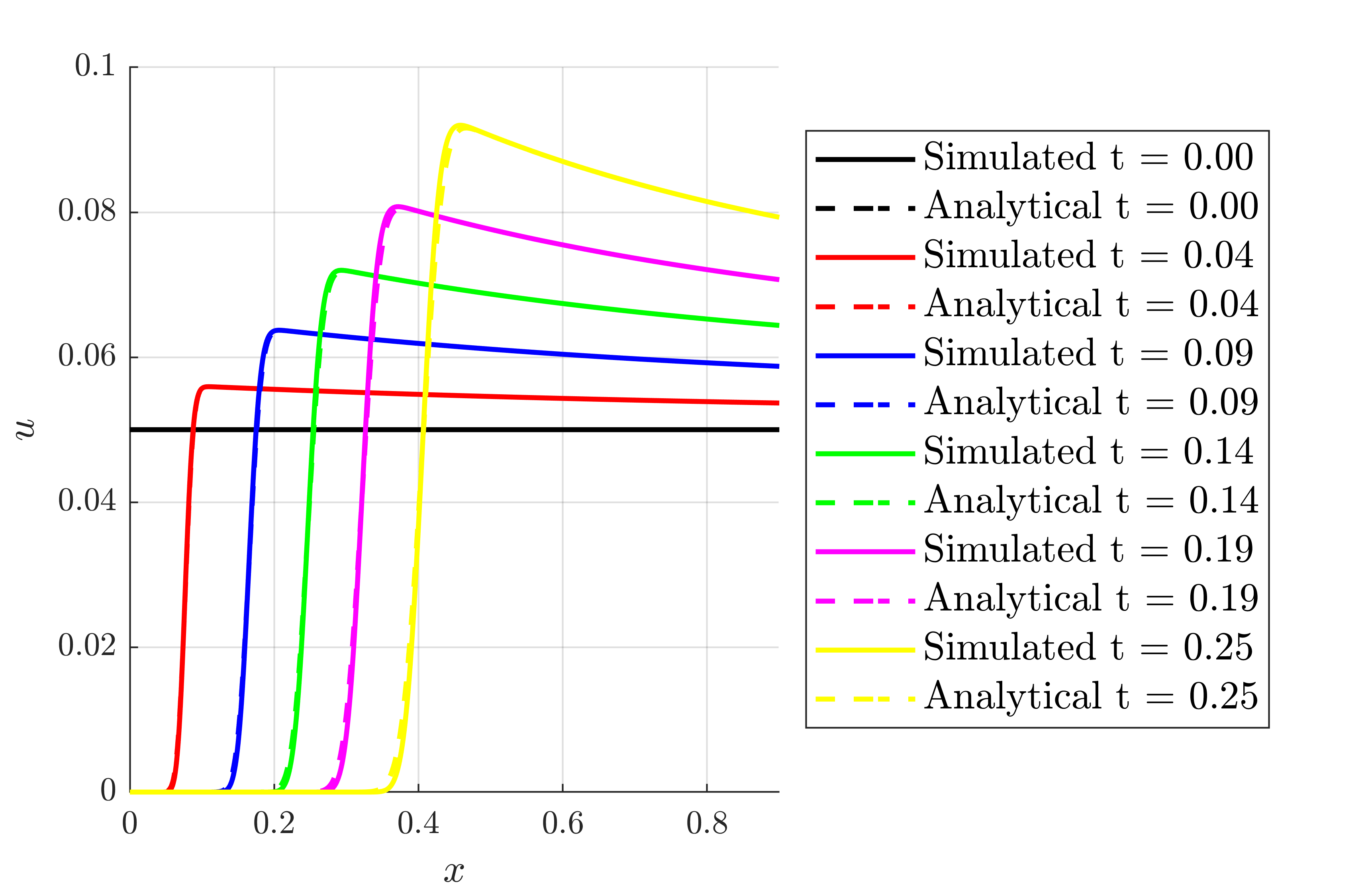}
    \caption[Comparison of numerical and analytical solutions for $\alpha$ exponential.]{\textbf{Comparison of numerical and analytical solutions for the exponential case.} The analytical and simulated profiles are compared at different times, considering $\beta(t,x)=1=g(t)$, $\alpha(t,x) = a\exp(-\lambda x)$ with $a = 2$,   $\lambda = 1 $ and the initial condition $u_0(x)=u_0=0.05 \ll 1$. In that case, we use $D = 1 \times 10^{-3}$ for computing the numerical solutions.}\label{fig:exp_comp}
\end{figure}

\subsubsection{Exponential chemotaxis}

We consider now a function of the type:
\begin{equation}
    \alpha(t,x) = \left(a\exp(-\lambda x)+b\exp(\lambda x)\right)g(t).
\end{equation}

Then, we have $f(x)=a\exp(-\lambda x)+b\exp(\lambda x)$, with Eq.~(\ref{eq:x_sol}) reducing to 
\begin{equation}
\arctan\left(\frac{b\exp(\lambda x)}{\sqrt{ab}}\right) - \arctan\left(\frac{b\exp(\lambda  s)}{\sqrt{ab}}\right) = \sqrt{ab}\lambda \int_0^\tau g(\eta)\,\mathrm{d}\eta  = \sqrt{ab}\cal{T}(\tau).
\end{equation}
Hence 
\begin{subequations} \label{eq:FG_exp}
\begin{align}
    F(t;s) &= \frac{1}{\lambda}\ln\left(\frac{\sqrt{ab} + a\tan\left(\sqrt{ab}\lambda \intgt \right)\mathrm{e}^{-\lambda s}}{\sqrt{ab}\mathrm{e}^{-\lambda s} -b\tan\left(\sqrt{ab}\lambda \intgt \right)}\right), \\
    G(t;x) &= \frac{1}{\lambda}\ln\left(\frac{\sqrt{ab} - a\tan\left(\sqrt{ab}\lambda \intgt \right)\mathrm{e}^{-\lambda x}}{\sqrt{ab}\mathrm{e}^{-\lambda x}  + b\tan\left(\sqrt{ab}\lambda \intgt \right)}\right). 
\end{align}
\end{subequations}
with 
the transition is located at
\begin{equation} \label{eq:xstar_exp}
    x^*(t) = \frac{1}{\lambda}\ln\left(\frac{\sqrt{ab} + a\tan\left(\sqrt{ab}\lambda \intgt \right)}{\sqrt{ab} -b\tan\left(\sqrt{ab}\lambda \intgt \right)}\right).
\end{equation}

In particular, if $b = 0$ this reduces to 
\begin{subequations} \label{eq:FG_exp_part}
\begin{align}
    F(t;s) &= \frac{1}{\lambda}\ln\left(a\lambda \intgt  + \mathrm{e}^{\lambda s}\right), \\
    G(t;x) &= \frac{1}{\lambda}\ln\left(-a\lambda \intgt  + \mathrm{e}^{\lambda x}\right), 
\end{align}
\end{subequations}
and
\begin{equation} \label{eq:xstar_exp_part}
    x^*(t) = \frac{1}{\lambda}\ln\left(1+\lambda a \intgt \right).
\end{equation}
    
Fig. \ref{fig:exp_comp} shows the comparison between the numerical results and the analytical solution with $\alpha(t,x) = 2\exp(-x)$ and $\beta(t,x) = 1=g(t)$, $D=1\times 10^{-3}$. Furthermore, in this regime  on  
neglecting $\mathcal{O}(u_0^2)$ corrections,  we have 
$$ u(x,t) = u_0\mathrm{e}^t \left. \frac{ f(s)}{f(x)} \right|_{s=G(t;x)} = \frac{u_0 \mathrm{e}^t}{1-a\lambda t \mathrm{e}^{-\lambda x}}, $$ 
noting that $x\geq x_*(t)$ on this characteristic, so that 
$$ 1-a \lambda t \mathrm{e}^{-\lambda x} \geq 1- a \lambda t \mathrm{e}^{-\lambda x_*(t)} = 1- \frac{a\lambda t}{1+a\lambda t}= \frac 1 {1+a\lambda t} >0, 
$$
in turn demonstrating that $u(x,t)$ does not possess a singularity with exponential decay in $x$ for $t$ fixed away from the transition to a good approximation, as may be also observed in Fig. \ref{fig:exp_comp}.

\clearpage
\section{Applications to microfluidic experiments} \label{sec:app}

\subsection{The general model} \label{tgm} 

We study a broad class of problems that are related to the evolution of a cell culture in microfluidic devices under a chemotactic agent, such as oxygen, when the concentration of the agent can be computed or measured, as schematically  represented in Fig. \ref{fig:scheme}. Hence we proceed to 
 consider the following model of generic cell culture evolution in microfluidic devices:

\begin{subequations} \label{eq:governing}
\begin{align}
    \frac{\partial C_n}{\partial T} &=  \frac{\partial}{\partial X}\left(D\frac{\partial C_n}{\partial X} - \chi C_n\frac{\partial B}{\partial X}\right) + \alpha_n M_g(B)C_n\left(1-\frac{C_n}{c_\mathrm{sat}}\right) - \alpha_{nd}M_d(B)C_n, \label{eq:governing_a} \\
    \frac{\partial C_d}{\partial T} &=  \alpha_{nd}M_d(B)C_n, \label{eq:governing_b} \\
    \frac{\partial B}{\partial T} &=  \frac{\partial}{\partial X}\left(D_B\frac{\partial B}{\partial X}\right) - \alpha_B  W(B,C_n). \label{eq:governing_c}
\end{align}
\end{subequations}
where $C_n$ and $C_d$ are respectively the alive and dead cell concentrations and $B$ is a chemotactic agent, $M_g$ and $M_d$ are nonlinear dimensionless corrections accounting for the effect of the chemo-attractant on cell growth and death. In addition $W(B,C_n)$ is the nonlinear dimensionless correction of the chemo-attractant consumption by all cells in the microdevice, noting that additional cells other than those of direct interest may be present, such as the study of metastatic tumour cells of interest migrating within a stromal cell population for instance. We also assume that dead cell concentration is  sufficiently low that it does not  compromise either live cell proliferation or migration. Note that if the chemotaxis agent is a cell nutrient, such as glucose or oxygen, $W(B,C_n) \neq 0$, whereas for unconsumed biochemical signals, $W(B,C_n)=0$.  The boundary conditions considered are
\begin{subequations} \label{eq:boundary}
\begin{align}
    \left.\frac{\partial f_n}{\partial x}\right|_{x=0} &= 0, \\
    \left.\frac{\partial f_n}{\partial x}\right|_{x=L} &= 0, \\
    B(x=0,t) &= B^L(t), \\
    B(x=L,t) &= B^R(t),
\end{align}
\end{subequations}
with $B^L(t)$ and $B^R(t)$, the chemotactic agent concentration at the left and right channel, and $f_n = D\frac{\partial C_n}{\partial X} - \chi C_n\frac{\partial B}{\partial X}$ the alive cell flow.  As the dead cell population is slave to the other variables in this model, we neglect it henceforth  which is equivalent to taking $\alpha_{nd} \simeq 0$.

Therefore, the full model here analysed is
\begin{subequations} \label{eq:igoverning}
\begin{align}
    \frac{\partial C_n}{\partial T} &=  \frac{\partial}{\partial X}\left(D\frac{\partial C_n}{\partial X} - \chi C_n\frac{\partial  B}{\partial X}\right)  +  \alpha_n M_g(B) C_n \left(1-\frac{C_n}{c_{\mathrm{sat}}}\right), \label{eq:igoverning_a} \\
    \frac{\partial B}{\partial T} &=  \frac{\partial}{\partial X}\left(D_B\frac{\partial B}{\partial X}\right) - \alpha_B  W(B,C_n). \label{eq:igoverning_b}
\end{align}
\end{subequations}

In order to evaluate the relevance of the different phenomena, we define the dimensionless variables:

\begin{subequations} \label{eq:dimensionless_variables}
\begin{align}
    C_n &=  c_\mathrm{sat}u, \label{eq:dimensionless_variables_u} \\
    B &=  B^*v, \label{eq:dimensionless_variables_v} \\
    X &=  Lx, \label{eq:dimensionless_variables_x} \\
    T &=  \frac{t}{\alpha_n}, \label{eq:dimensionless_variables_t}
\end{align}
\end{subequations}
Hence Eqs. (\ref{eq:igoverning}) become
\begin{subequations} \label{eq:governing_d}
\begin{align}
    \frac{\partial u}{\partial t} &=  \frac{\partial}{\partial x}\left(\Pi_1\frac{\partial u}{\partial x} - \Pi_2 u\frac{\partial v}{\partial x}\right) + m(v)u\left(1-u\right), \label{eq:governing_a_d} \\
    \frac{\partial v}{\partial t} &=  \frac{\partial}{\partial x}\left(\Pi_3\frac{\partial v}{\partial x}\right) + \Pi_4 w(v,u), \label{eq:governing_b_d}
\end{align}
\end{subequations}
where
\begin{subequations} \label{eq:Pi}
\begin{align}
    \Pi_1 &=  \frac{D}{\alpha_n L^2}, \label{eq:Pi_1} \\
    \Pi_2 &=  \frac{\chi\mathrm{O_2^*}}{\alpha_n L^2}, \label{eq:Pi_2} \\
    \Pi_3 &=  \frac{D_B}{\alpha_n L^2}, \label{eq:Pi_3} \\
    \Pi_4 &=  \frac{\alpha_B c_{\mathrm{sat}}}{ \alpha_n B^*}, \label{eq:Pi_4} \\
    m(v)  &=  M_g(B^*v), \label{eq:r} \\
    w(v,u)  &=  W(B^*v,c_{\mathrm{sat}}u). \label{eq:s} \\
\end{align}
\end{subequations}
The associated boundary conditions are 
\begin{subequations} \label{eq:boundary_d}
\begin{align}
    \left.\frac{\partial f_n}{\partial x}\right|_{x=0} &= 0, \label{eq:boundary_d1} \\
    \left.\frac{\partial f_n}{\partial x}\right|_{x=1} &= 0, \label{eq:boundary_d2} \\
    v(x=0,t) &= \psi_1(t), \label{eq:boundary_d3}\\
    v(x=1,t) &= \psi_2(t), \label{eq:boundary_d4}
\end{align}
\end{subequations}
where now $f_n = \Pi_1\frac{\partial u}{\partial x} - \Pi_2 u\frac{\partial v}{\partial x}$ and $\psi_1(t),~\psi_2(t)$ are prescribed functions of time corresponding to the level of nutrient or chemoattractant fixed to be at the channel edges.  In what follows, we will use $u_t$ and $u_x$ as an abbreviation for $\frac{\partial u}{\partial t}$ and $\frac{\partial u}{\partial x}$.

In particular the governing PDEs described in Eqs.~(\ref{eq:governing_d}) may be reformulated as 
\begin{subequations} \label{eq:governing_dimensionless_1}
\begin{align}
    u_t &=  \Pi_1 u_{xx} -\Pi_2 \left(v_x u\right)_x + m(v)u\left(1-u\right), \label{eq:governing_dimensionless_1a} \\
    v_t &=  \Pi_3 v_{xx} - \Pi_4 w(v,u). \label{eq:governing_dimensionless_1b}
\end{align}
\end{subequations}

\subsubsection{The weak consumption limit} \label{twcl}

First we consider the case where the chemoattractant is not a nutrient and therefore is not consumed by cells. In that case, we can set $$w(v,u) = 0.$$ For all the problems in which it is possible to assume $\Pi_1 \ll 1$,  whereby  random cellular motility is negligible compared to directed chemotaxis and with $\Pi_3 \gg 1$,  so that the  chemoattractant diffusion is large relative to cellular diffusion Eqs.~(\ref{eq:governing_dimensionless_1}) reduce to 
\begin{subequations}
\begin{align}
    u_t  + k \left(v_x u\right)_x &=  m(v)u\left(1-u\right), \label{eq:governing_dimensionless_2a} \\
    v_{xx} &= 0. \label{eq:governing_dimensionless_2b}
\end{align}
\end{subequations}
where $k = \Pi_2$.

However, care is required in considering the boundary conditions for $u$ and the initial conditions for $v$ in this reduced model, due to the loss of the second spatial derivative of $u$ and the first temporal derivative of $v$. In particular, we cannot satisfy all the boundary conditions for $u$; instead we have boundary layers. We have seen in the examples these occur  at internal transitions and at the right of the domain (see Fig \ref{fig:osc_comp}). Thus for the simplified system the boundary condition is enforced at $x=0$ for $u$ with the simplification of the flux to $f_n=-kv_x u$, as diffusion is treated as negligible. However,  enforcing the boundary condition at $x=1$ will require the consideration of a boundary layer that is not resolved in the simplified model as it is complicated, but does not further insight into cell migratory and chemotaxis. The two boundary conditions for $v$, as given by Eqs. \ref{eq:boundary_d3} and \ref{eq:boundary_d4}  are inherited and applied at both boundaries. We now consider the initial conditions for $v$, which cannot be satisfied. Instead there is an analogous temporal boundary layer for early time while initial transients relax, though such transients persist for such a short time that they are not of interest, and thus not resolved, here. The justification of the neglect of these fast transients is  further detailed in Appendix \ref{appft}, where it is demonstrated that the solution of Eqn (\ref{eq:governing_dimensionless_2b}) corresponds to the leading order composite solution in a temporal boundary layer analysis that exploits $\Pi_3\gg1$ .

Proceeding, we set $v(x=0,t) = \psi_1(t)$ and $v(x=1,t) = \psi_2(t)$, Eq.~(\ref{eq:governing_dimensionless_2b}) is immediately integrated to \begin{eqnarray}\label{vout1}
v(t,x)=\left(\psi_2(t) - \psi_1(t)\right)x + \psi_1(t). 
\end{eqnarray} 
We recover therefore Eq. (\ref{eq:eq_governing_simp}) with $$\alpha(t,x) = k(\psi_2(t) - \psi_1(t)),~~~~~~~~~\beta(t,x)  = m(\left(\psi_2(t) - \psi_1(t)\right)x + \psi_1(t)),$$ that is a special case of the linear problem $\alpha(t,x) = a(t)x+b(t)$, with $a(t) = 0$ and $b(t) = k(\psi_2(t) - \psi_1(t))$, so the different functions needed in order to compute the solution are:
    
    \begin{subequations}
    \begin{align}
        F(t;s) &=  s + k\int_0^t\Delta \psi(\eta)\, \mathrm{d}\eta, \label{eq:ex1_1} \\
        G(t;x) &= x -  k\int_0^t\Delta \psi(\eta)\, \mathrm{d}\eta, \label{eq:ex1_2} \\
        x^*(t) &=  k\int_0^t\Delta \psi(\eta)\, \mathrm{d}\eta = x-G(t;x),  \label{eq:ex1_3}
    \end{align}
    \end{subequations}
    where we have defined $\Delta \psi (t) = \psi_2(t) - \psi_1(t)$. Also, the expression of the cell profile far from the transition is
    \begin{equation}
        u(x,t) =  \frac{u_0( G(t;x))\exp\left(\int_0^t K(\eta, F(\eta;G(t;x))) \, \mathrm{d}\eta\right)}{1+u_0( G(t;x))\left(\exp\left(\int_0^t K(\eta, F(\eta;G(t;x)))\, \mathrm{d}\eta\right) - 1\right)} \label{eq:ex1_4}
    \end{equation}
    where
    \begin{equation}
        K(\eta,X) = m\left((\psi_2(\eta) - \psi_1(\eta))X + \psi_1(\eta)\right).
    \end{equation}

\begin{figure}[!t]
\centering
\begin{subfigure}{\textwidth}
  \centering
  \includegraphics[clip=true,trim=0pt 0pt 0pt 0pt,width=0.6\textwidth]{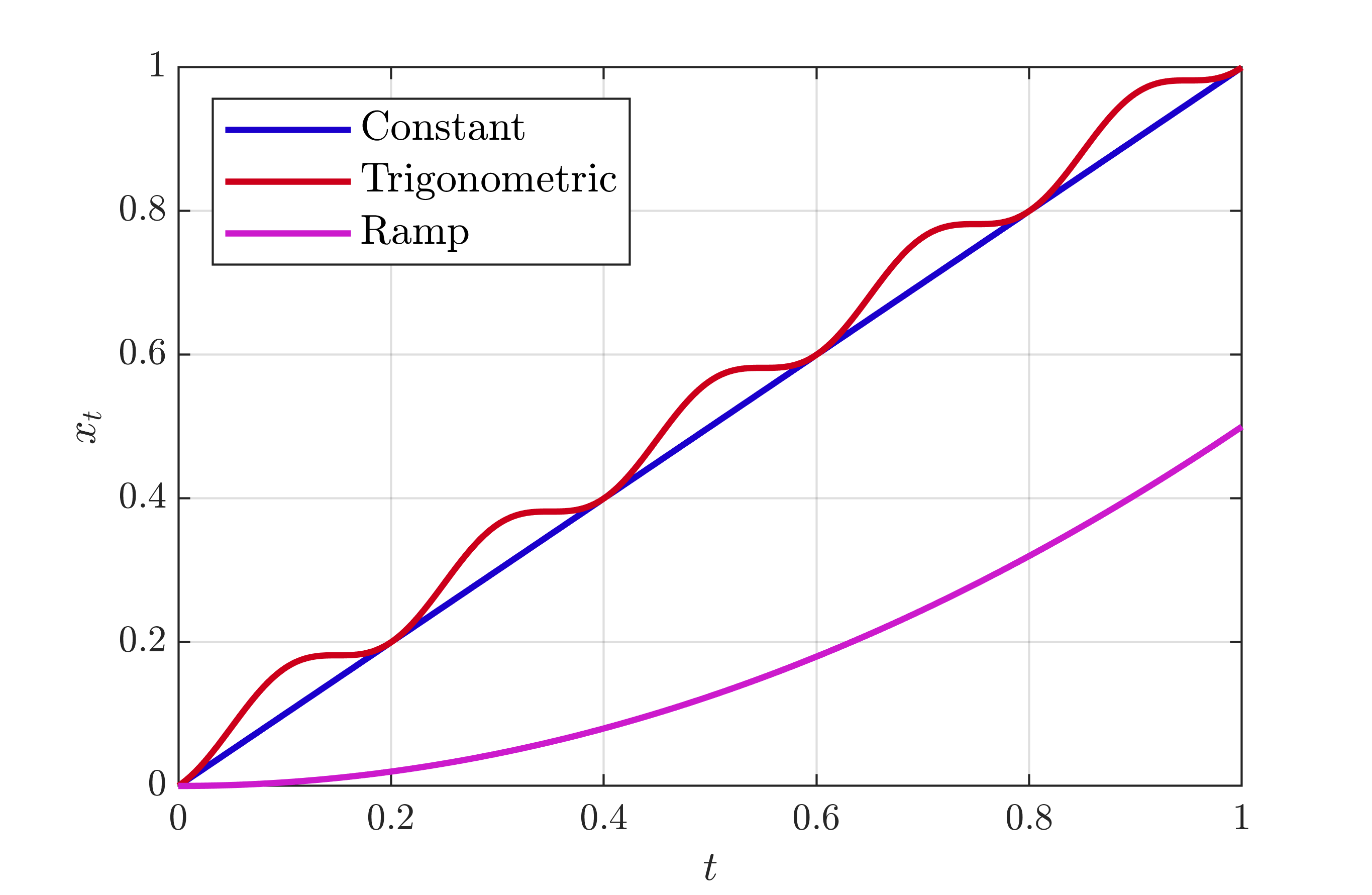}
  \caption{Value of the transition coordinate $x^* = x^*(t)$.}
  \label{fig:solution_1a}
\end{subfigure}
\begin{subfigure}{\textwidth}
  \centering
  \includegraphics[clip=true,trim=0pt 0pt 0pt 0pt,width=0.6\textwidth]{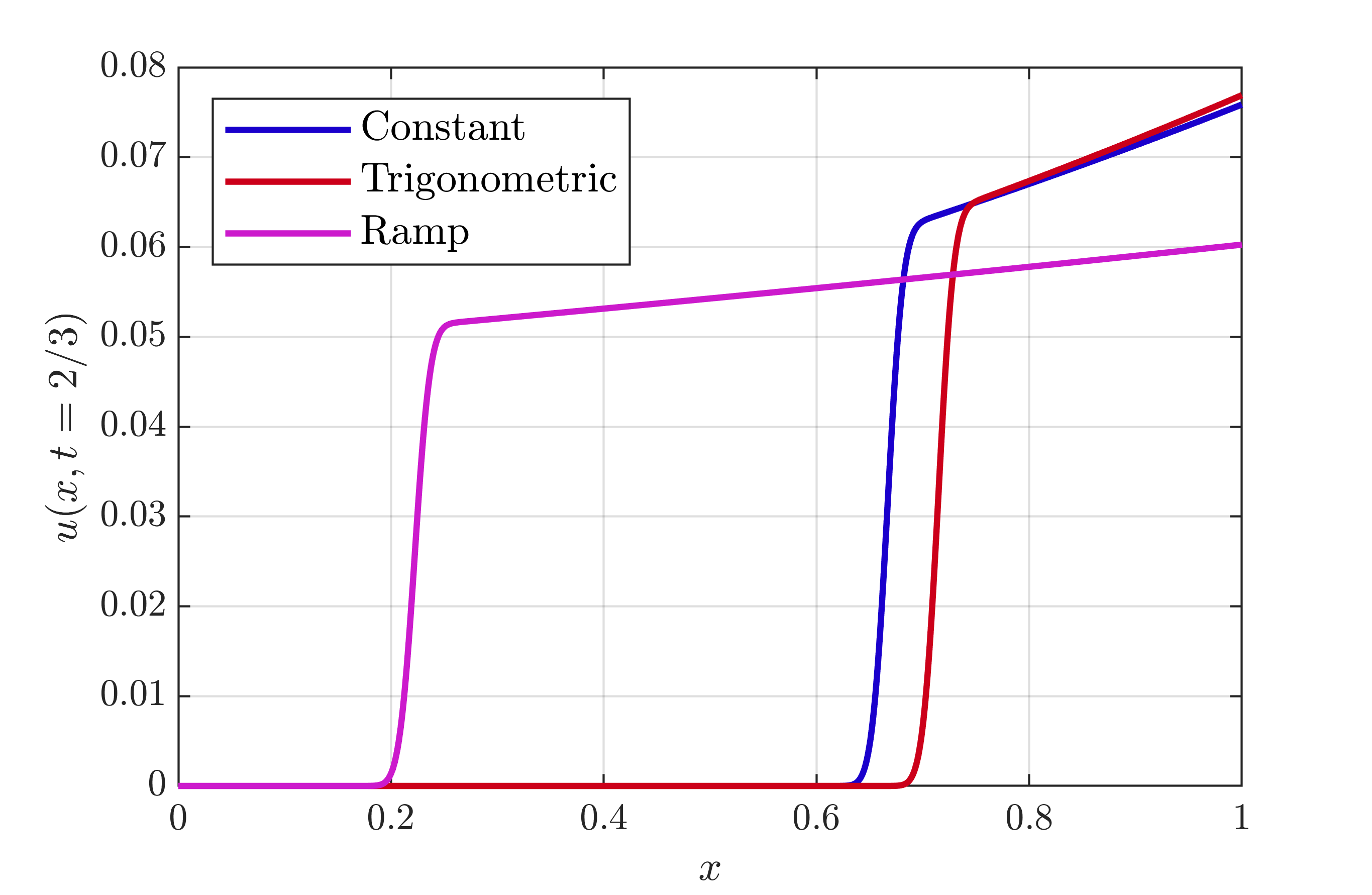}
  \caption{Cell profile at $t=2/3$.}
  \label{fig:solution_1b}
\end{subfigure}
\caption[Solution for the case with no consumption]{\textbf{Approximate analytical solution for the case with no consumption.} We assume $k=1$, $m(v) = v$ and the transition region is recreated with $D = 1 \times 10^{-4}$. The initial condition is set again to $u_0(x) = u_0 = 0.05$. Also, we consider three different shapes for $\psi_2$, with $\psi_1=0$. A constant oxygen level $\psi_2(t) = 1$, a trigonometric oscillatory stimulus $\psi_2(t) = 1 + \cos(\omega t)$ with $\omega = 10\pi $, and an increasing ramp stimulus $\psi_2(t) = t$.}
\label{fig:solution_1}
\end{figure}

The evolution of the transition coordinate $x^*(t)$ and the dimensionless cell profile for different times are shown in Fig. \ref{fig:solution_1} for $k=1$, $\psi_1(t) = 0$ and different external stimuli $\psi_2(t)$. In particular, let us consider the case with $m(v) = m_0 v$, and with $\psi_1(t)=0$ and use of the change of variable $X=F(\eta;s)$, whereby on a characteristic 
$$\dfrac{\mathrm{d}X}{\mathrm{d}\eta } = \alpha(\eta,X)=k\psi_2(\eta)$$
with $F(t;s)=x$ and $F(0,s)=s$. This reveals 
\begin{eqnarray*} \left. \int_0^t K(\eta,F(\eta,s))\, \mathrm{d}\eta \right|_{s=G(t;x)} &=&  \left. \int_s^x \frac 1 {k\psi_2(\eta)}  m_0 \psi_2(\eta) X \, \mathrm{d} X \right|_{s=G(t;x)}= \frac{m_0}{2k} (x-G(t;x))(x+G(t;x)) \\ &=&  \frac{m_0}{2k} x^*(t)(2x-x^*(t)). 
\end{eqnarray*}
Combined with Eq. \eqref{eq:ex1_4}, this gives a simple expression for $u(x,t)$.  For instance 
$$ u(x,t) \approx u_0 \exp\left( \frac{m_0}{k} x^*(t)(x-x^*(t)/2) \right) +\mathcal{O}(u_0^2),$$
for $u_0\ll 1$ constant;  in this case, we have an increasing function at fixed $t$ to the right of the transition given $x^*(t)>0$ and this is essentially  linear for $x^*(t)\ll 1$, as observed  in Fig. \ref{fig:solution_1}.

\subsubsection{Cellular consumption of chemoattractant}\label{ccoc} 

A more interesting case is when the chemoattractant is a nutrient and therefore it is consumed by cells. In that case, $w(v,u) \neq 0$. 

With $w_{\mathrm{pc}}(v)$ denoting the non-dimensional uptake of nutrient {\it per cell}, we take $w_{\mathrm{pc}}(v)$ to be monotonic increasing with 
\begin{subequations}
\begin{align}
    \lim_{v\rightarrow 0} w_{\mathrm{pc}}(v) = 0, \\
    \lim_{v\rightarrow +\infty} w_{\mathrm{pc}}(v) = 1, 
\end{align}
\end{subequations}
where the final limit is without loss of generality, with the overall scale of uptake governed by $\Pi_4$. 

For instance, with Michaelis-Menten kinetics we take  \cite{cornish2013origins}:
\begin{equation}
    w_{\mathrm{pc}}(v) = \frac{v}{v+k_m},
\end{equation}
or, more in general, Hill-Langmuir equation for modelling the consumption kinetics \cite{hill1910possible} 
\begin{equation}
    w_{\mathrm{pc}}(v) = \frac{v^n}{v^n+k^n_H},
\end{equation}

\begin{figure}[!htbp]
\centering
\begin{subfigure}{.49\textwidth}
  \centering
  \includegraphics[clip=true,trim=30pt 170pt 100pt 80pt,width=\textwidth]{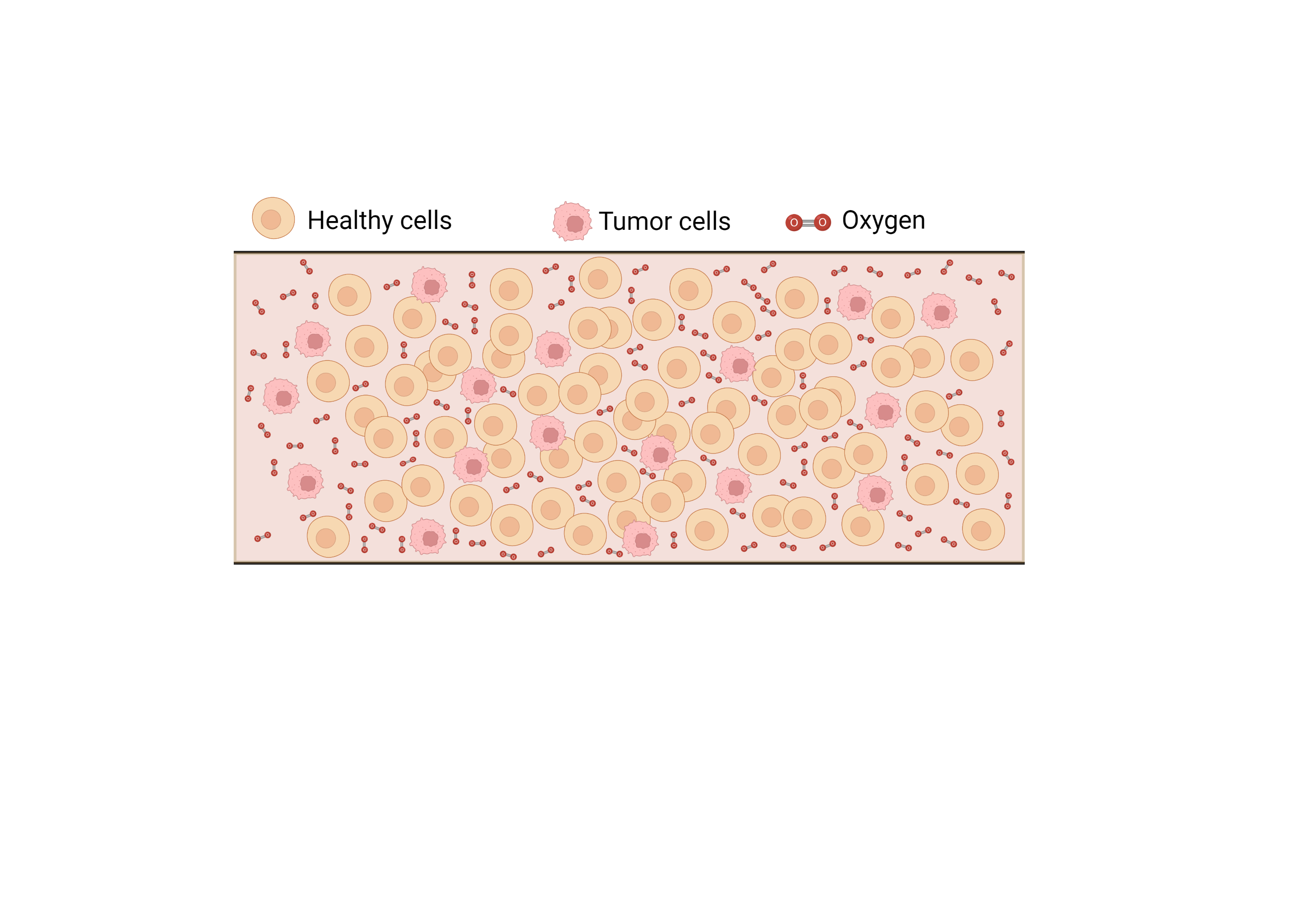}
  \caption{High oxygen levels.}
  \label{fig:illustration_a}
\end{subfigure}
\begin{subfigure}{.49\textwidth}
  \centering
  \includegraphics[clip=true,trim=30pt 170pt 100pt 80pt,width=\textwidth]{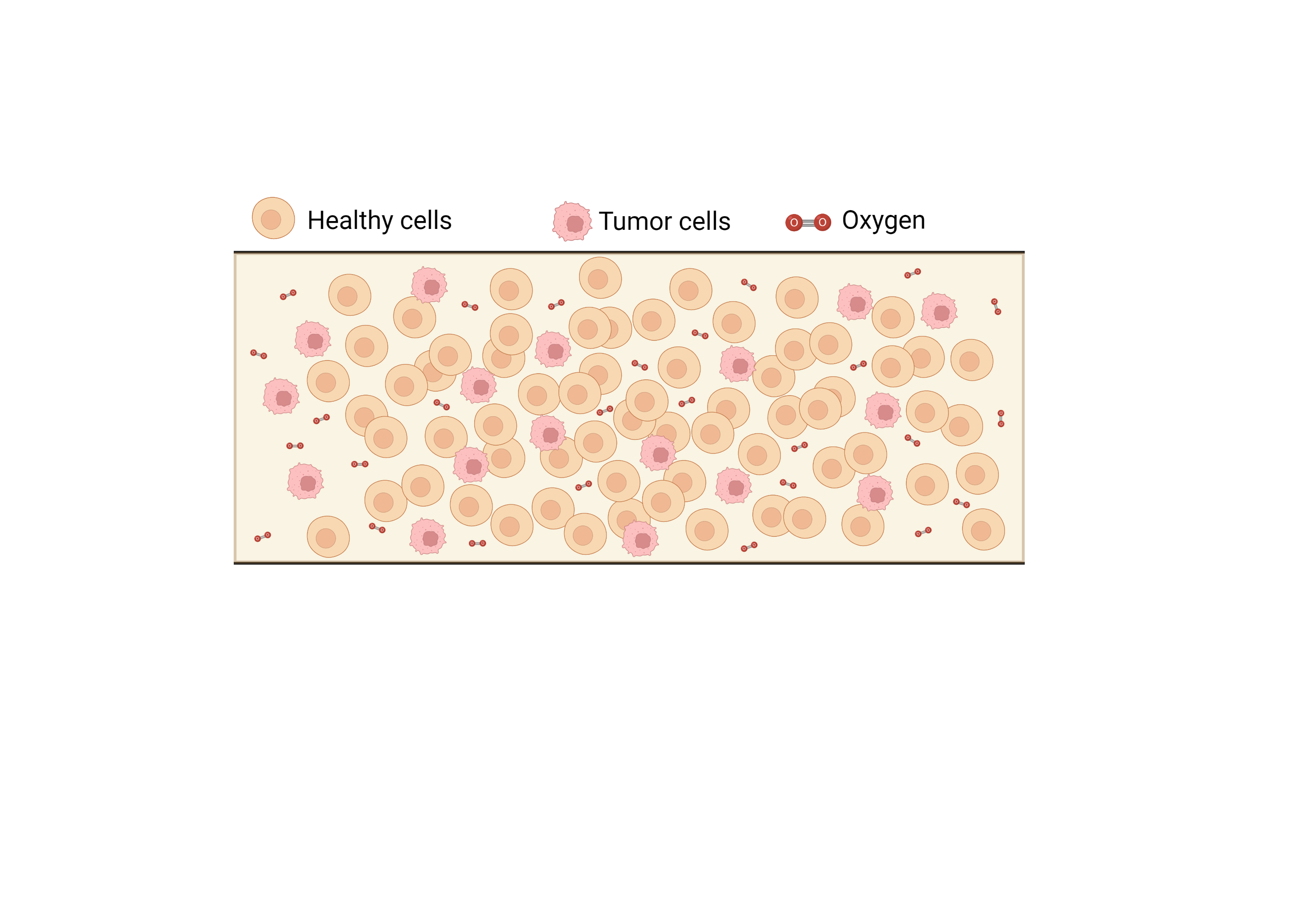}
  \caption{Low oxygen levels.}
  \label{fig:illustration_b}
\end{subfigure}
\caption[Cell culture model recreating cancer cells in an oxygenated ambient.]{\textbf{Cell culture model recreating cancer cells in an oxygenated ambient.} For illustration purposes, the oxygen is considered as the chemoattractant  nutrient of the heterogeneous growth and death. The different local and global oxygen gradients appearing in the chamber may be approximated by approximating Eq.~(\ref{eq:governing_dimensionless_1b}) when considering specific situations, such as high oxygen levels or low oxygen levels. Created with BioRender.com.}
\label{fig:illustration}
\end{figure}

In any of the aforementioned cases, there are numerous potential scenarios:
\begin{enumerate}
    \item There are no other cells at the microfluidic device besides the cell culture of our interest and we are at the low cell regime, $\Pi_4 u/\Pi_3 \ll 1$. In that case, after rapid initial transients   describing the diffusion relaxation of the nutrient, and too fast to be on the timescale of cellular motility,  Eq.~(\ref{eq:governing_dimensionless_1b}) becomes $v_{xx} = 0$, and the discussion is analogous to the case without the consumption term.
  
    \item There are other non-migrating cells within the microfluidic device in addition to the migrating cells, for instance if we are considering  a metastasis model, with the other cells at  constant concentration and in excess of the tumour cells. If we additionally have high nutrient concentrations the situation is illustrated in Fig. \ref{fig:illustration_a} using oxygen as an example of nutrient. In that case, $w_{\mathrm{pc}}(v) \sim 1$ and 
    have 
    $$ w(u,v) = K,$$ 
    where  $K$ is the (dimensionless) total amount of cells, essentially constant as the non-tumour cells are in excess. Then, with the definition     $$ \lambda = \Pi_4 K/\Pi_3,$$   Eq. (\ref{eq:governing_dimensionless_1b}) 
       becomes 
       \begin{eqnarray} \label{vint1}  v_t = \Pi_3 v_{xx} - \Pi_3 \lambda . \end{eqnarray}
       In addition, we assume $\Pi_3^{-1}\ll \lambda \ll \Pi_3,$ so that  $\lambda$ may be treated as order one on using asymptotic methods based on  on the leading order of approximations based on $\Pi_3\gg 1$. Then the above further  reduces to
       \begin{eqnarray} \label{vfin1}  v_{xx} =   \lambda ,  \end{eqnarray}
     noting, as above,  that fast initial transients are not of interest, with   further justification of Eqn \eqref{vfin1} in Appendix \ref{appft} via a boundary layer analysis.

    \item There are other cells within the microfluidic device at constant concentration and in excess relative to the tumour cells, together with low nutrient concentrations. The situation is also illustrated in Fig. \ref{fig:illustration_b}. Assuming the Michaelis-Menten model, $w_\mathrm{pc}(v) \sim v/k_m$ so that 
    $w(u,v) = Kv/k_m$ where $K$ is the effectively  constant non-dimensional total cell density and  Eq.~(\ref{eq:governing_dimensionless_1b}) 
      reduces to 
     \begin{eqnarray} \label{vint2}  v_t = \Pi_3 v_{xx} - \Pi_3 \lambda v,  \end{eqnarray}
     where now $\lambda =   \Pi_4 K  /(k_m\Pi_3) $. As above,  this reduces if  $\Pi_3^{-1}\ll \lambda \ll \Pi_3,$ which we assume in order to yield \begin{eqnarray} \label{vfin2}  v_{xx} =   \lambda v,  \end{eqnarray}
    once more noting  fast initial transients are not of interest,  with additional justification of Eqn \eqref{vfin2} presented in Appendix \ref{appft}.

\end{enumerate}

\paragraph{The presence of other cells in excess and high chemoattractant or oxygen levels.} 

With the assumptions and motivations as above in the derivation of Eq. \eqref{vfin1} we now have 
\begin{subequations}
\begin{align}
    u_t  + k \left(v_x u\right)_x &=  m(v)u\left(1-u\right), \label{eq:governing_dimensionless_3a} \\
    v_{xx} &= \lambda. \label{eq:governing_dimensionless_3b}
\end{align}
\end{subequations}
where the boundary conditions are again the ones given by Eqs.~(\ref{eq:boundary_d}), except for the fact that the flux is given by $f_n=-kv_x u$ as diffusion is neglected and the cell boundary condition at $x=1$, with its associated boundary layer, is no longer considered as discussed in detail previously.  

If we prescribe $v(x=0,t) = \psi_1(t)$ and $v(x=1,t) = \psi_2(t)$, Eq.~(\ref{eq:governing_dimensionless_3b})  integrates to 
\begin{eqnarray} \label{vout2} v(x,t)=\frac{1}{2}\lambda x^2 + \left(\psi_2(t) - \psi_1(t) - \frac{1}{2}\lambda\right)x + \psi_1(t).
\end{eqnarray} 
We have therefore $$\alpha(t,x) = k\lambda x + k\left(\psi_2(t) - \psi_1(t)-\frac{1}{2}\lambda\right),
$$
and $$ \beta(t,x) = m\left(v(t,x)\right),$$ for cell growth.  Once more, this is a special case of $\alpha(t,x) = a(t)x+b(t)$ with $a(t) = k\lambda$ and $b(t) = k\left(\psi_2(t) - \psi_1(t)-\frac{1}{2}\lambda\right)$, so the different functions needed in order to compute the solution are:

    \begin{subequations}
    \begin{align}
        F(t;s) &= s\mathrm{e}^{k\lambda t} + \frac{1}{2} \left(1-\mathrm{e}^{k\lambda t}\right) + k\int_0^t\Delta \psi(\eta)\mathrm{e}^{-k\lambda \eta} \, \mathrm{d}\eta, \label{eq:ex2_1} \\
        G(t;x) &= x\mathrm{e}^{-k\lambda t} - \frac{1}{2} \left(\mathrm{e}^{-k\lambda t}-1\right) - k\mathrm{e}^{-k\lambda t}\int_0^t\Delta \psi(\eta)\mathrm{e}^{-k\lambda \eta} \, \mathrm{d}\eta, \label{eq:ex2_2} \\
        x^*(t) &= \frac{1}{2} \left(1-\mathrm{e}^{k\lambda t}\right) + k\int_0^t\Delta \psi(\eta)\mathrm{e}^{-k\lambda \eta}\, \mathrm{d}\eta. \label{eq:ex2_3}
    \end{align}
    \end{subequations}

The expression of the cell profile may be computed using Eq.~(\ref{eq:u_sol}), which gives 

\begin{equation}
        u(x,t) = \frac{u_0(G(t;x))\exp(-k\lambda t)\exp\left(\int_0^t K(\eta, F(\eta;G(t;x))) \, \mathrm{d}\eta \right)}{1 + u_0(G(t;x))\int_0^t K(\eta, F(\eta;G(t;x)))\exp(-k\lambda \eta) \exp\left(\int_0^\eta K(\xi, F(\xi;G(t;x))) \, \mathrm{d}\xi \right)\, \mathrm{d}\eta} ,
\end{equation}
where now
\begin{equation}
        K(\eta,X) = m\left(\frac{1}{2}X^2 + (\psi_2(\eta) - \psi_1(\eta) - \frac{1}{2}\lambda)X + \psi_1(\eta)\right).
\end{equation}

\begin{figure}[!t]
\centering
\begin{subfigure}{\textwidth}
  \centering
  \includegraphics[clip=true,trim=0pt 0pt 0pt 0pt,width=0.6\textwidth]{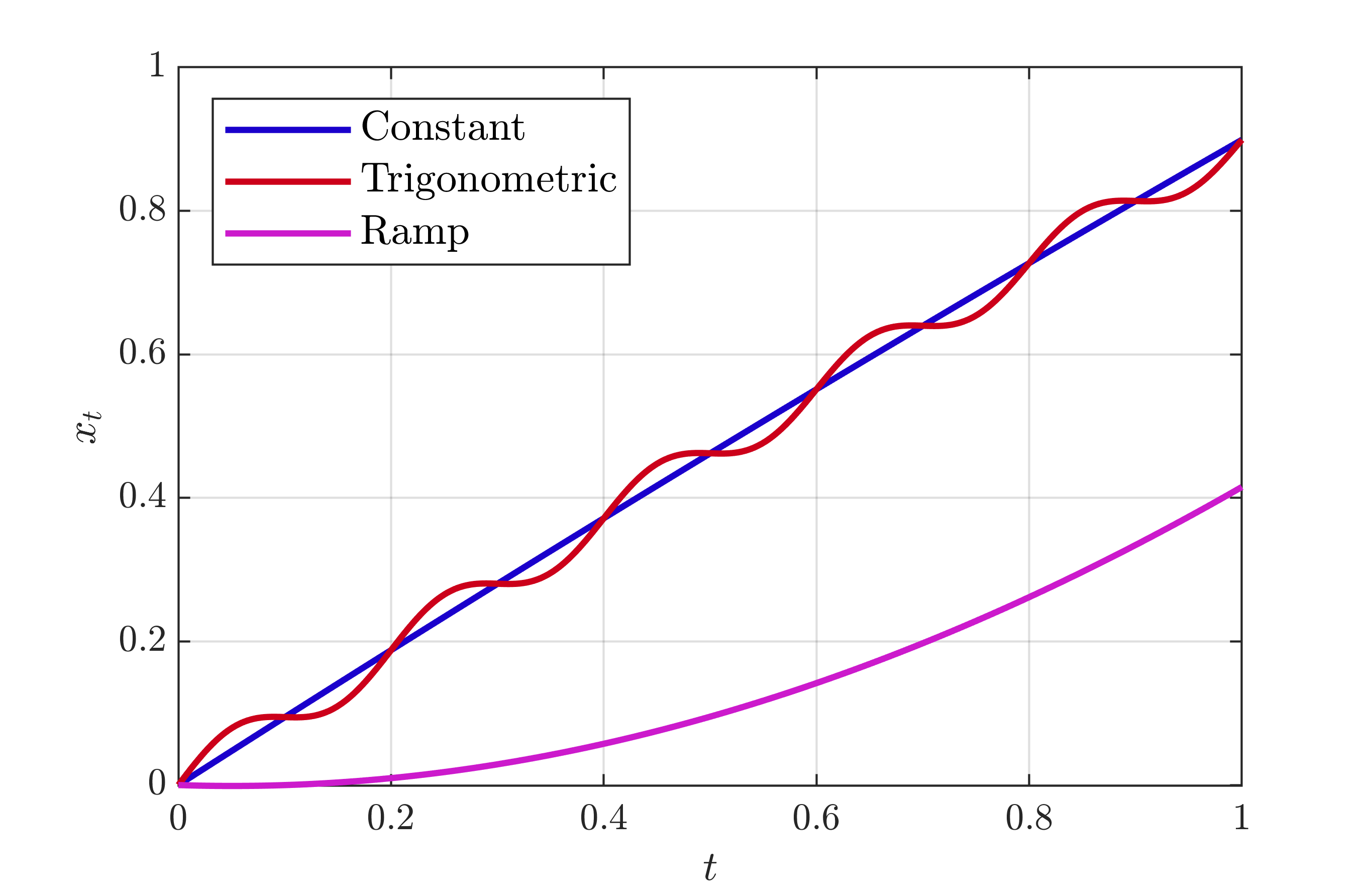}
  \caption{Value of the transition coordinate $x^* = x^*(t)$.}
  \label{fig:solution_2a}
\end{subfigure}
\begin{subfigure}{\textwidth}
  \centering
  \includegraphics[clip=true,trim=0pt 0pt 0pt 0pt,width=0.6\textwidth]{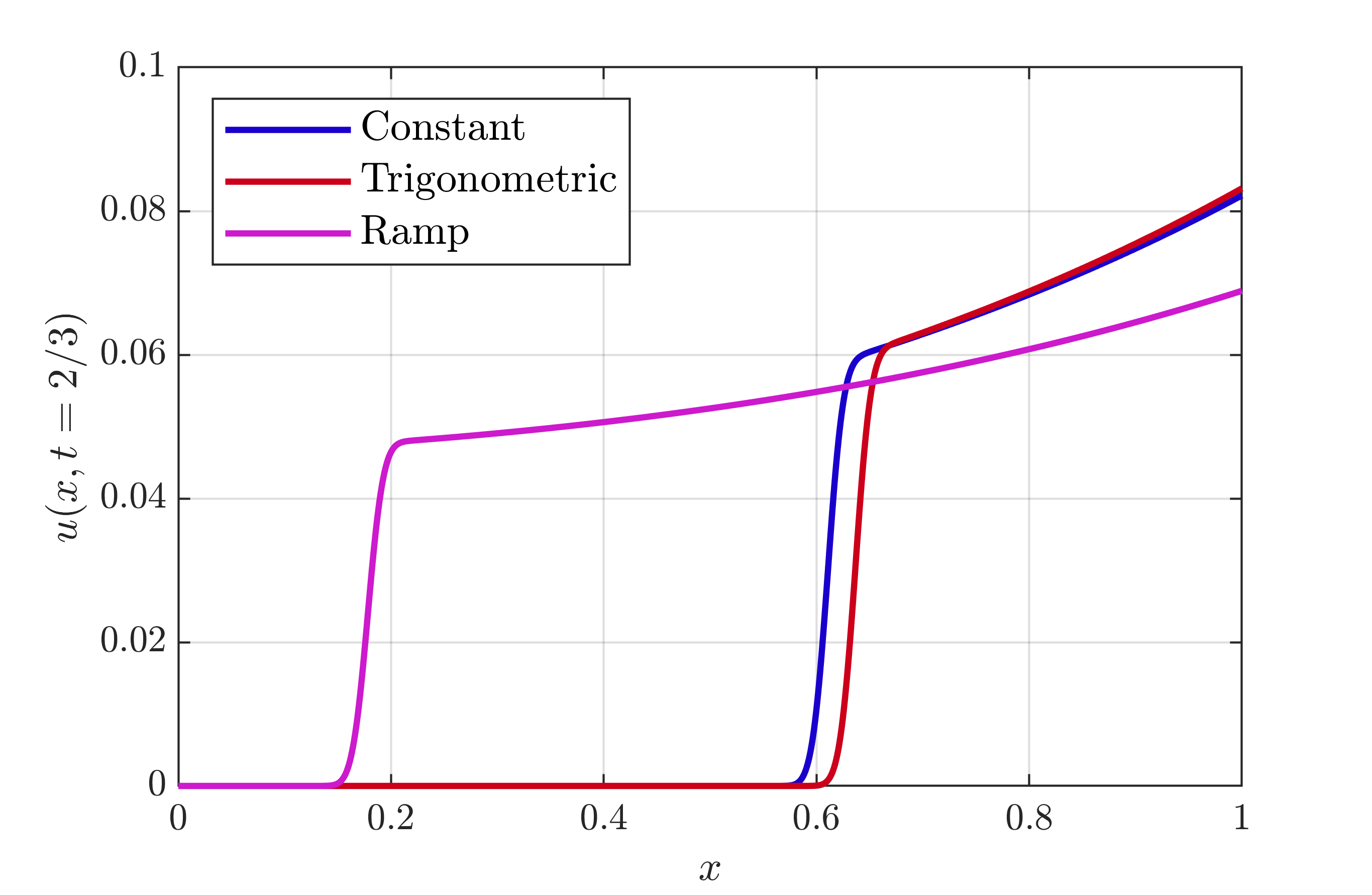}
  \caption{Cell profile at $t=2/3$.}
  \label{fig:solution_2b}
\end{subfigure}
\caption[Solution for the case with chemoattractant consumption and high chemoattractant levels]{\textbf{Solution for the case with chemoattractant consumption and high chemoattractant levels.} Approximate analytical solution for $k=1$,   $m(v)=v$,  $\lambda = 0.1$ and the transition region is recreated with $D = 1 \times 10^{-4}$. The initial condition is set again to $u_0(x) = u_0 = 0.05$. We consider three different shapes for $\psi_2$. A constant oxygen level $\psi_2(t) = 1$, an oscillatory stimulus $ \psi_2(t) = 1 + \cos(\omega t)$ with $\omega = 10\pi$, and an increasing stimulus $\psi_2(t) = t$.}
\label{fig:solution_2}
\end{figure}

The evolution of the transition coordinate $x^*(t)$ and the dimensionless cell profile for different times are shown in Fig. \ref{fig:solution_2} for $k=1$, $\lambda=0.1$, $m(v) = m_0v$ with $m_0 = 1$, $u_0 = 0.05$ and a external stimulus given by $\psi_1(t) = 0$ and different shapes for $\psi_2(t)$. 

While a fast oscillatory stimulus with $\psi_2(t)=1+\cos(10 \pi t)$  plotted in Fig. \ref{fig:solution_2} does not allow a ready approximation for cell density to the right of the transition region, we can consider the cases of $\psi_2(t)=1$ or $\psi_2(t)=t$ that are also considered in this Figure. In particular, we can use  elementary but extensive manipulation to determine and approximate 
$$ \left. \int_0^t K(\eta, F(\eta,s)) \, \mathrm{d}\eta \right|_{s=G(t;x)},$$ to deduce that 
$$ u_0 \approx u_0 \mathrm{e}^{-\lambda t} \exp\left( \frac {1-\mathrm{e}^{-2\lambda t}}{4\lambda} (x-x^*(t))^2 + \frac 1 3 h_1(t,\lambda)(x-x^*(t)) +  h_2(t,\lambda) \right) + \mathcal{O}(u_0^2),
$$ in these cases. 
In particular, $h_1(t,\lambda)$ is of the form $$h_1(t,\lambda) = \left( 3t + \frac 3 2 t^2\right)-\left(\frac 3 2 + \frac 9 4 t +   t^2 \right) \lambda t +
\mathcal{O}(\lambda^2t^2),$$
for $\psi_2(t)=1$ and $$h_1(t,\lambda) = \left(\frac{3}{2}t^2 + \frac{1}{2}t^3\right) - \left(\frac{3}{2}+ \frac{3}{4}t + \frac{1}{2}t^2+\frac{3}{8}t^3\right)\lambda t + \mathcal{O}(\lambda^2t^2),$$
for $\psi_2(t)=t$.

For  $h_2(t,\lambda)$, we have  with $\psi_2(t)=1$ that  $$ h_2(t, \lambda) =   \frac 1 2 t^2 \left( 1+\frac 1 3 t \right) - t \left(  \frac 1 2 + \frac 1 3 t + \frac 1 8 t^2  \right)\lambda t 
+ \mathcal{O}(\lambda^2t^2),$$
while, in contrast, for $\psi_2(t)=t$ we note that $$h_2(t,\lambda) =  \frac{t^4}{40} \left(5 + t\right)  - t^2\left(\frac{1}{4}+\frac{1}{16}t+\frac{1}{15}t^2+\frac{1}{36}t^3\right)\lambda t + \mathcal{O}(\lambda^2t^2).$$ 

For both cases note that at fixed time the cell concentration to the right of the transition is essentially the exponential of a quadratic in $x-x^*(t)$, provided that $\lambda$ is small enough. 

\paragraph{Other cells and low chemoattractant or oxygen levels.} With analogous reasoning, Eqs. (\ref{eq:governing_dimensionless_1}) now reduce to 
\begin{subequations}
\begin{align}
    u_t  + k \left(v_x u\right)_x &=  m(v)u\left(1-u\right), \label{eq:governing_dimensionless_4a} \\
    v_{xx} &= \lambda v. \label{eq:governing_dimensionless_4b}
\end{align}
\end{subequations}

The general solution to Eq.~(\ref{eq:governing_dimensionless_4b}) generates  a chemotactic function $\alpha(t,x)$ that is neither linear in $x$ nor separable, with Eq. \ref{eq:cF2b_general} equivalent to a Riccati differential equation in  $\exp(x)$, for which general solutions are not known in terms of standard functions. Even though the general case is not tractable in terms of constructing solutions, two particular important configurations do allow progress, namely the gradient configuration ($v(x=0,t)=0$ and $v(x=1,t)= \psi(t)$) and the symmetric configuration (($v(x=0,t)=v(x=1,t)=\psi(t)$).

The gradient configuration is certainly the most interesting configuration. For that case, Eq.~(\ref{eq:governing_dimensionless_4b}) integrates  to 
\begin{eqnarray} \label{vout3} v(t,x)=\psi(t)\sinh(\sqrt{\lambda} x),
\end{eqnarray}  whereby  
$$\alpha(t,x) = k\sqrt{\lambda}\cosh(\sqrt{\lambda} x)\psi(t), ~~~~~~~~~ \beta(t,x) = m_0(v(t,x)).$$  Hence we have separability, with  $f(x) = k\sqrt{\lambda}\cosh(\sqrt{\lambda} x)$ and $g(t) = \psi(t)$. With 
$$ {\cal S}(t) = \tan(k\lambda \intgt /2), ~~~~~~~ \intgt = \int_0^t g(\eta)\mathrm {d}\eta$$ 
we in turn have
\begin{subequations}
\begin{align}
    F(t;s) &= \frac 1{\sqrt{\lambda}}\ln\left(\frac{1+\mathrm{e}^{-s\sqrt{\lambda}}{\cal S}(t)}{\mathrm{e}^{-s\sqrt{\lambda}}-{\cal S}(t)}\right), \label{eq:ex3a_1} \\
    G(t;x) &= \frac 1{\sqrt{\lambda}}\ln\left(\frac{1-\mathrm{e}^{- x\sqrt{\lambda}}{\cal S}(t)}{\mathrm{e}^{- x\sqrt{\lambda}}+{\cal S}(t)}\right),  \label{eq:ex3a_2} \\
    x^*(t) &= \frac 1{\sqrt{\lambda}}\ln\left(\frac{1+{\cal S}(t)}{1-{\cal S}(t)}\right). \label{eq:ex3a_3}
\end{align}
\end{subequations}

As with the previous cases the expression of the cell profile may be computed using Eq.~(\ref{eq:u_sol}), whereby 

\begin{equation}
        u(x,t) =  \frac{u_0(G(t;x))\exp\left(\int_0^t K(\eta, F(\eta;G(t;x))) \, \mathrm{d}\eta \right)}{1 +  u_0(G(t;x))\int_0^t K(\eta,F(\eta;G(t;x))) \exp\left(\int_0^\eta K(\xi, F(\xi;G(t;x))) \, \mathrm{d}\xi  \right)\, \mathrm{d}\eta }, \label{equfin}
\end{equation}
where $K(\eta,X) = m(v) - kv_{xx}$, that is
\begin{equation}
        K(\eta,X) = m\left(\psi(\eta)\sinh(\sqrt{\lambda} X)\right) - k \lambda\psi(\eta)\sinh(\sqrt{\lambda} X),
\end{equation}
and, with $m(v) = m_0 v$ and $u_0(x)=u_0 \ll 1,$ constant, use of the change of variables $X=F(\eta,s)$ reduces  Eq. \eqref{equfin} to 
$$ u(x,t) = u_0 \left(\frac{\cosh(\sqrt{\lambda}G(t;x)) }{\cosh(\sqrt{\lambda} x)}   \right)^{\frac{m_0-k\lambda}{k\lambda} } +\mathcal{O}(u_0^2). 
$$

\begin{figure}[!t]
\centering
\begin{subfigure}{\textwidth}
  \centering
  \includegraphics[clip=true,trim=0pt 0pt 0pt 0pt,width=0.6\textwidth]{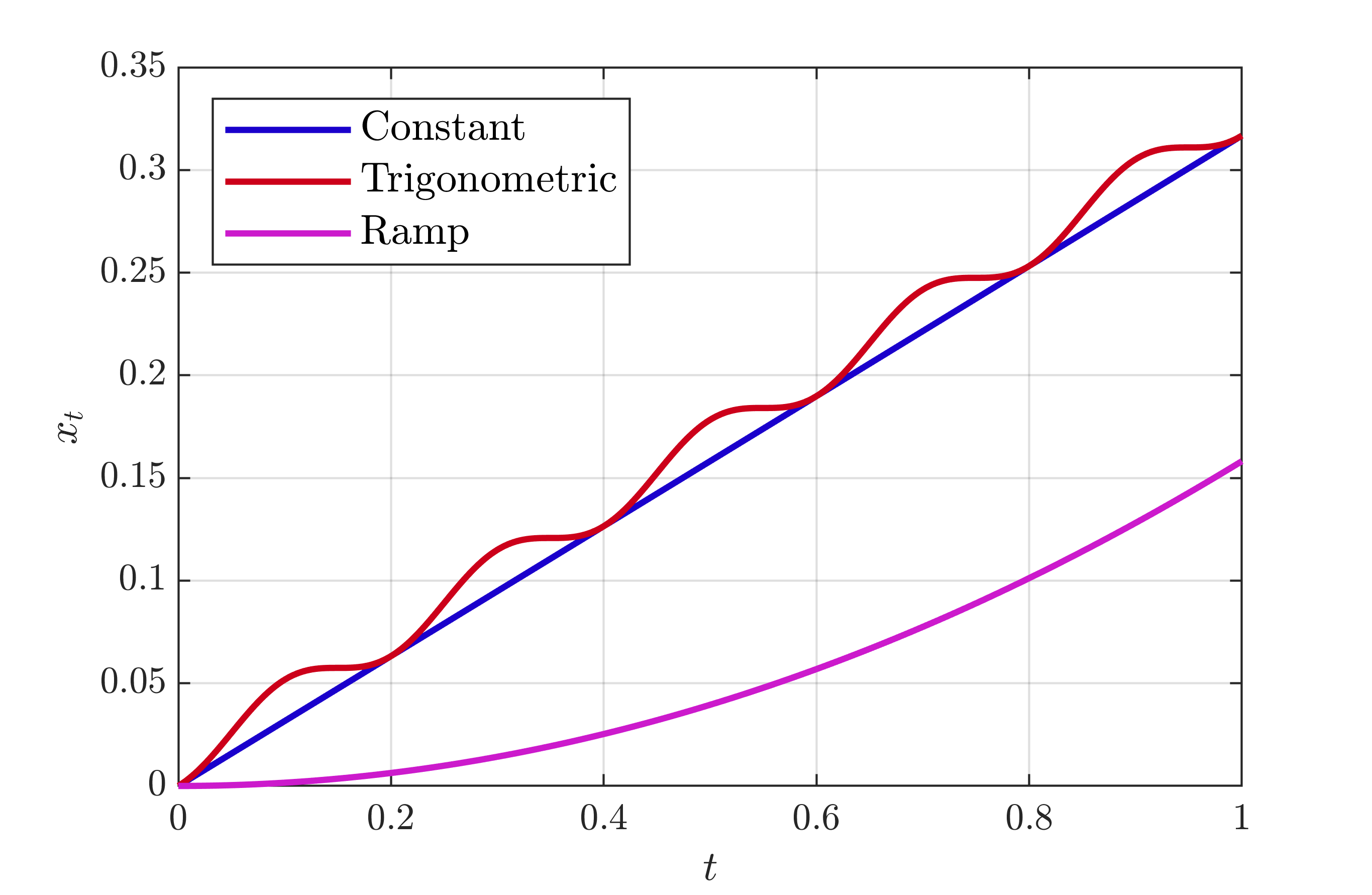}
  \caption{Value of the transition coordinate $x^* = x^*(t)$.}
  \label{fig:solution_3a}
\end{subfigure}
\begin{subfigure}{\textwidth}
  \centering
  \includegraphics[clip=true,trim=0pt 0pt 0pt 0pt,width=0.6\textwidth]{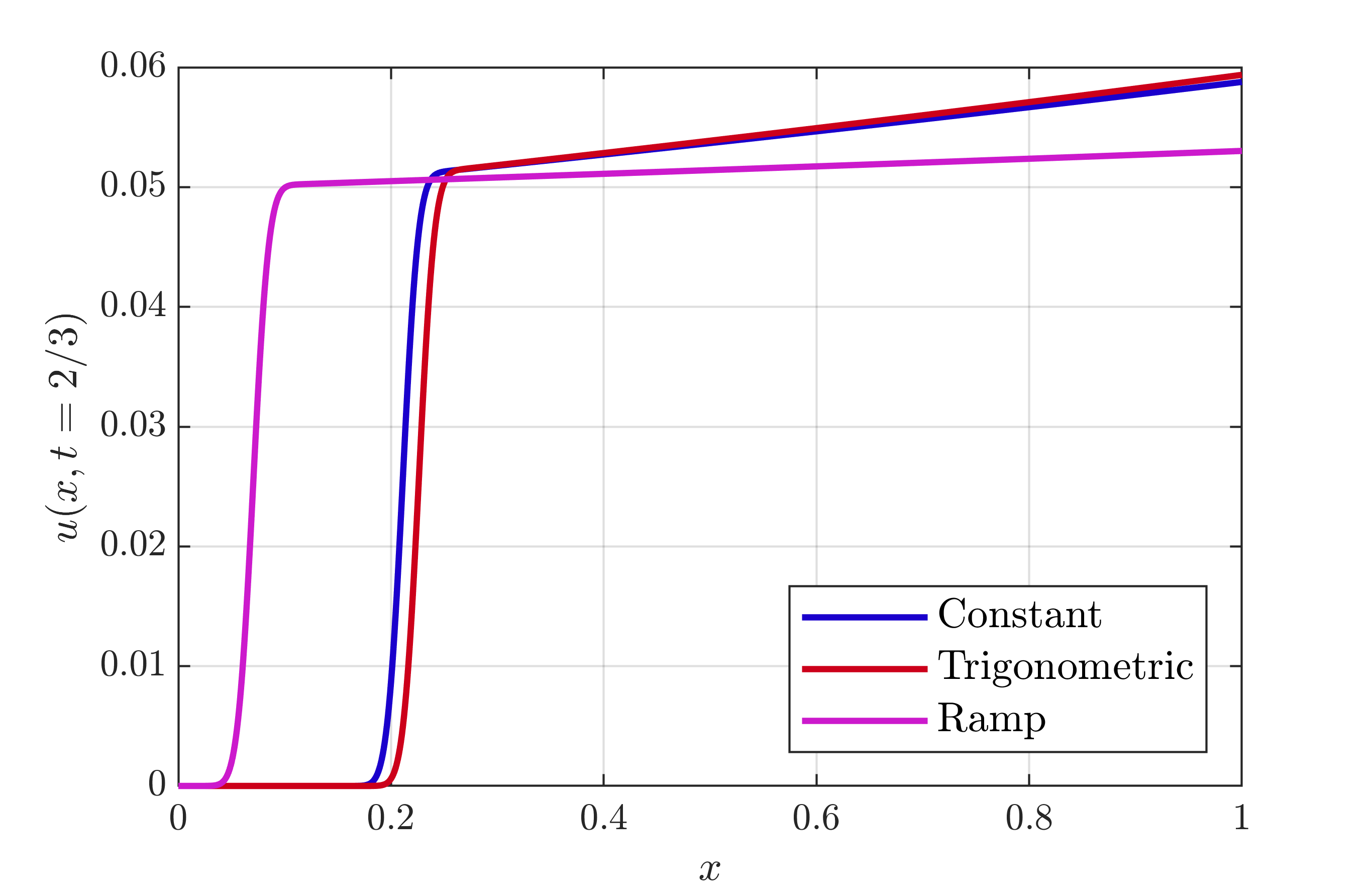}
  \caption{Cell profile at $t=2/3$.}
  \label{fig:solution_3b}
\end{subfigure}
\caption[Solution for the case with chemoattractant consumption and low chemoattractant levels]{\textbf{Solution for the case with chemoattractant consumption and low chemoattractant levels.} Approximate analytical solution for $k=1$, $m(v) = v$, $\lambda = 0.1$ and the transition region is recreated with $D = 1 \times 10^{-4}$. The initial condition is set again to $u_0(x) = u_0 = 0.05$. We consider three different shapes for $\psi_1$. A constant oxygen level $\psi_1(t) = 1$, an oscillatory stimulus $\psi_1(t) = 1 + \cos(\omega t)$ with $\omega = 10\pi $, and an increasing stimulus $\psi_1(t) = t$.}
\label{fig:solution_3}
\end{figure}

The evolution of the transition coordinate $x^*(t)$ and the dimensionless cell profile for different times are shown in Fig. \ref{fig:solution_3} for $m_0=1$, $k=1$, $\lambda = 0.1$, $ u_0(x)=u_0=0.05,$ constant. Here,  the cellular density to the right of transition further simplifies to 
$$ u(x,t) =  u_0 \left( \frac{1+{\cal S}(t)^2}{(1-{\cal S}(t)^2) +2 {\cal S}(t)\sinh\left(\frac{x}{\sqrt{10}}\right) }\right)^{9} +\mathcal{O}(u_0^2), ~~~~~ {\cal S}(t) = \tan(\intgt / 20),$$ which entails for small time the spatial variation of the cellular density is that of a hyperbolic sine. Finally,  note from the expression for $x^*(t)$ and the fact $x^*(t) \leq 1$ for the transition zone to be within the domain, we have the bound ${\cal T} \leq \tanh(1/2) < 1 $ so that the denominator in the above expression is positive and there is no singularity.

Even if less interesting from the experimental point of view, we can also obtain an expression for a symmetric configuration by considering half of the domain and applying Neumann boundary conditions for $x=0$ so Eq.~(\ref{eq:governing_dimensionless_4b}) is integrated to 
\begin{eqnarray} \label{vout4}v(t,x)=\psi(t)\cosh(\sqrt{\lambda} x).  \end{eqnarray} 
Again we are under a separable case with $\alpha(t,x) = k\sqrt{\lambda}\sinh(\sqrt{\lambda} x)\psi(t)$ and $\beta(t,x) = 1$, so $f(x) = k\sqrt{\lambda}\sinh(\sqrt{\lambda} x)$ and $g(t) = \psi(t)$, so the different functions needed in order to compute the solution are:

\begin{subequations}
\begin{align}
F(t;s) &= \frac 2 {\sqrt{\lambda}} \arctan\left(\mathrm{e}^{k\lambda\intgt}\tanh\left(\frac{s\sqrt{\lambda}}{2}\right)\right), \label{eq:ex3b_1} \\
G(t;x) &= \frac 2 {\sqrt{\lambda}}
\arctan\left(\mathrm{e}^{-k\lambda\intgt}\tanh\left(\frac{x\sqrt{\lambda}}{2}\right)\right), \label{eq:ex3b_2} \\
x^*(t) &= 0. \label{eq:ex3b_3} 
\end{align}
\end{subequations}

Now in Eq. (\ref{equfin}) we shall use
\begin{equation}
       K(\eta,X)= m\left(\psi(\eta)\cosh(\sqrt{\lambda} X\right) - k \lambda\psi(\eta)\cosh(\sqrt{\lambda} X).
\end{equation}
and, with $m(v) = m_0 v$ and  $u_0(x)=u_0 \ll 1,$ constant, Eq. (\eqref{equfin}) reduces to 
$$ u(x,t) = u_0 \left(\frac{\sinh(\sqrt{\lambda}x) }{\sinh(\sqrt{\lambda} G(t;x))}   \right)^{\frac{m_0-k\lambda}{k\lambda} } +\mathcal{O}(u_0^2). 
$$ 

% \textcolor{red}{In this last equation, the fraction should be replaced by its inverse, am I right? EAG ... I make this way around, have sent calculation by email}

\clearpage
\section{Discussion}

There is an increasing use of \emph{in vitro} investigations in the exploration of cellular motility, for instance with microdevice studies exploring tumour cell dynamics in response to oxygen gradients, as illustrated by GBM studies \cite{ayuso2016development,ayuso2017glioblastoma,ayensa2020mathematical}. In turn this has motivated the main theme of this paper, namely modelling cell migration chemotaxis in heterogeneous environments, which applies for general chemoattractants, not just oxygen. In particular, the governing equations for cellular motility that have been considered are of the form 
\begin{equation}
    \frac{\partial u}{\partial t} = D \frac{\partial^2 u}{\partial x^2} - \frac{\partial }{\partial x}\left(\alpha(t,x)u_x\right) + \beta(t,x)u(1-u),
\end{equation}
and supplemented by zero flux boundary conditions and suitable, typically  constant, initial conditions. 
Furthermore, the general spatiotemporal heterogeneity in the chemotactic function $\alpha(t,x)$ and the growth function $\beta(t,x)$ emerges from the chemoattractant gradients that the cells respond to, which may be manipulated extensively in microdevice experiments. However, the heterogeneity also entails that while cell migration dynamics towards high concentrations of chemoattractant is anticipated, the dynamics will not simply be that of a translationally invariant travelling wave and their associated analytical simplicity.

Hence, we have considered a  framework capable of considering spatiotemporal chemotactic gradients and the resulting wavefront and cell density dynamics for  cellular migrations in the presence of spatial and temporal heterogeneity. 
In particular, the fact cell spreading in the absence of chemotaxis is generally much slower than in its presence, the non-dimensional cellular diffusion coefficient is very small, so that typically  $D\ll 1$, which we assume and exploit in this study.  With this, we have that away from boundary layers near sources, which were located at the domain edge in the examples considered, and away from cell wavefronts with sharp transitions,  an advection equation for the cellular density $u$ emerges with $D\sim 0$. 
Nonetheless, this advection equation is still complex, entailing that a constraint on constructing solutions for the cell density behaviour within  the analytical framework presented here is that the chemotactic function $\alpha(t,x)$ must either be linear in $x$ or separable with respect to space and time. However, numerous cases are consistent with these constraints, as illustrated by the range of examples considered  in Section \ref{spcoi}, together with  the examples that emerge from the consideration of cellular dynamics within microdevices in Section \ref{sec:app}. 

With these constraints, we have illustrated how the method of characteristics can be used  to construct the cell density solutions away from boundary layers and transition regions, together with the use of boundary layer methods to construct a uniformly continuous approximate solution that accommodates the transition in the wavefront of the cells. In turn this provides an analytical characterisation of the movement of cellular wavefronts and the concentrations of cells  within and either side of the wavefront. However we do not resolve boundary layers near oxygen sources, which has been the righthand  boundary in the examples considered, for the presented framework given the limited insight the boundary layer will provide for the overall cell behaviour. 

Even with the restriction of linearity or separability of the chemotactic functions $\alpha(t,x)$, there is extensive freedom in the choices of $\alpha(t,x)$ and $\beta(t,x)$. Hence we first considered the predictions of the model for cell behaviour in exemplar test cases. In particular, we have explored and documented the behaviour of the cell density where the chemotactic function varies  linearly in space and non-trivially in time, as well as quadratically in space and finally  with exponential spatial behaviour. 
In addition to analytical investigations, we have also verified that the analytical solutions faithfully reproduce the behaviour of direct numerical solutions of the model, as may be observed extensively in Figures \ref{fig:lin_comp}-\ref{fig:exp_comp}. 
Furthermore the exemplar solutions are subsequently used to study models reduced from a full coupling between the cellular density and chemoattractant concentration with a microdevice setting, as documented in Section \ref{sec:app}. 
 Once more the resulting dynamics are analytically investigated and documented, with rationally based  analytical approximation for the location of the transition region and the cellular density presented. 

Such  analytical solutions are useful in numerous ways. For example, they provide oversight and insight for the system dynamics across parameter space. To illustrate, we consider the wavefront being driven by chemoattractant at the right hand boundary, so that $\psi_1(t)=0, v(t,1)=\psi_2(t)$. Then, with $k$ a measure of the strength of chemotaxis,   differentiating $x^*(t)$ from Eq. \eqref{eq:ex2_3} for the speed for the transition front gives 
 $$ \frac{\mathrm{d}x^*}{\mathrm{d}t}= -\frac 1 2 k\lambda \mathrm{e}^{k\lambda t} + k \psi_2(t)\mathrm{e}^{-k\lambda t}$$ 
 for the model of section \ref{ccoc}, with high chemoattractant or oxygen levels. Hence, in this case, increasing the chemoattractant/oxygen uptake, as representing by increasing $\lambda$, always slows the propagation of a rightmoving wave. In contrast, for the model with low chemoattractant or oxygen levels of  section  \ref{ccoc} differentiating $x^*(t)$ from Eq. \eqref{eq:ex3a_3} we obtain $$ \frac{\mathrm{d}x^*}{\mathrm{d}t}= \frac{\sqrt{\lambda}}{\mathrm{cos}\left(k\lambda \mathcal{T}\right)}\psi_2(t)$$ therefore revealing that increasing the chemoattractant/oxygen uptake always speeds up the propagation of rightmoving wave, illustrating the general deductions that may be made from the presented analytical solutions. 

As well as providing analytical characterisation of the systems behaviour the analytical approximations provide a means to very rapidly compute cell behaviour. Thus the approximate analytical solutions can support computationally intensive studies. A simple example would be a global sensitivity analyses over all parameters. In particular, rapid evaluation would be most useful for parameter estimation using experimental, often noisy, data especially if Bayesian techniques are considered as this require extensive simulation to provide posterior distributions, rather than optimisation techniques which only generate point estimates for parameters. A directly analogous example is Bayesian model selection, whereby the comparison of experimental data and model prediction is used to distinguish different model structures when these are not known a priori, such as different functional forms of $\alpha(t,x)$ or $\beta(t,x)$ representing different prospective growth and chemotactic responses for a given the tumour cell line in question. In particular both optimisation techniques and the Bayesian techniques are iterative, so that the use of a rational but rapid evaluation of  an approximation to an optimum in optimisation studies or a posterior for Bayesian techniques can in turn be used to restart the procedure with the full numerical model to further refine the results \cite{brown2022risking}.

In summary, we have developed a framework for the construction of analytical approximations for front behaviours and densities for cells undergoing chemotaxis in heterogeneous environments, as characterised by the chemotactic and growth functions $\alpha$ and $\beta$ respectively.  The resulting cellular waves of migration are not simple travelling waves due to the heterogeneity induced by the chemoattractant profiles. Nonetheless, numerous  features of the wavefronts have been extracted via rational approximation, such  as the location and speed of the propagating wave, together with the cellular density profile. These have been explored and validated for exemplars scenarios as well as investigated for models fundamentally motivated by experimental microdevices for observing cellular motility under a very wide range of conditions, even if  complete generality is not feasible for progress using constructive methods. Thus, the  solutions presented here not only provide insight into the behaviour of cellular motility under the influence of spatiotemporal chemotactic heterogeneity, but also highlight general behaviours and important mechanisms and  parameters, as well as providing a means of rapid evaluation in demanding computational studies, such as Bayesian parameter estimation and model selection. 

\section*{Aknowledgements}
The authors gratefully acknowledge the financial support  from the Spanish Ministry of Science and Innovation (MICINN), the State Research Agency (AEI), and FEDER, UE through the projects PID2019-106099RB-C44/AEI and PID2021-126051OB-C41, the Government of Aragon (DGA) and the ``Centro de Investigación Biomedical en Red en Bioingeniería, Biomateriales y Nanomedicina (CIBER-BBN)'', financed by the Instituto de Salud Carlos III with assistance from the European Regional Development Fund (FEDER).

\clearpage

\bibliographystyle{unsrt}  
\bibliography{references} %%% Remove comment to use the external .bib file (using bibtex).

\clearpage
\appendix
\section{Oscillatory solutions for homogeneous growth} \label{sec:a1}\label{sec11} 
 We consider  homogeneous growth, so that $\beta(t,x)$ is a constant independent of $x$ and $t$, together with  a temporally oscillating chemotactic response governed by $\alpha(t,x)$, so that Eq.~(\ref{eq:edo_osc_prev}) reduces to 
\begin{equation} \label{eq:a_ini}
    r'+\left(\beta-a\cos(\omega t)\right)r =\beta, \quad r(0) = r^*. 
\end{equation}
Hence the  cell concentration to the right of the transition region, but away from the boundary layer at $x=1$, is spatially constant with  a temporal oscillation, which we consider below   for four distinct parameter regimes. 

\subsection{Slow variations of the gradients, \texorpdfstring{$\bm{\omega \ll 1}$}.}
Let $T = \omega t$ and, without loss of generality, we consider the decomposition  $r(t) = y+z(\omega t)$. Then, Eq.~(\ref{eq:a_ini}) becomes 
\begin{equation} \label{eq:edo_truco_1}
    \omega z'(\omega t) + \frac{\mathrm{d}y}{\mathrm{d}t} + \left(\beta - a\cos(\omega t)\right)z(\omega t) + \left(\beta - a\cos(\omega t)\right)y = \beta. 
\end{equation}

Further, let $\omega z'(\omega t) + \left(\beta - a\cos(\omega t)\right)z(\omega t)=\beta$, so we have 
\begin{equation} \label{eq:edo_f}
    \omega z'(T) +\left(\beta - a\cos(T)\right)z(T) = \beta. 
\end{equation}

Since $\omega \ll 1$, we can approximate the solution of Eq.~(\ref{eq:edo_f}) by 
\begin{equation} \label{eq:sol_f}
    z(T) = \frac{\beta}{\beta -a\cos(T)} + \mathcal{O}(\omega). 
\end{equation}

Then, Eq.~(\ref{eq:edo_truco_1}) yields 
\begin{equation} \label{eq:edo_truco_2}
  \frac{\mathrm{d}y}{\mathrm{d}t}+ \left(\beta - a\cos(\omega t)\right)y = 0,  
\end{equation}
 and the initial condition is 
 \begin{equation}
     y^* = y(t=0) = r(t=0) - z(0) = r^* - \frac{\beta}{\beta - a}. 
 \end{equation}

As we have a slow modulation of the frequency/decay rate, we use  the  Wentzel-Kramers-Brillouin (WKB) method \cite{olver1997asymptotics}. In terms of $T = \omega t$, Eq.~(\ref{eq:edo_truco_2}) becomes 
\begin{equation} \label{eq:edo_truco_3}
   \omega \frac{\mathrm{d}y}{\mathrm{d}T} + \left(\beta - a\cos(T)\right)y = 0, \quad y(0) = y^*.  
\end{equation}

The WKB approximation is expressed here as 
\begin{equation}\label{eq:WKBJ}
    y = p\exp\left(\frac{\phi(T)}{\omega}\right)J(T), \quad J(T) = J_0 + \omega J_1 + \mathcal{O}(\omega^2).
\end{equation}

Substituting Eq.~(\ref{eq:WKBJ}) into Eq.~(\ref{eq:edo_truco_2})  we obtain 
\begin{equation}
    p\exp\left(\frac{\phi(T)}{\omega}\right)\left[\omega\left(\frac{\dot{\phi}}{\omega}J + \dot{J}\right)+(\beta - a\cos T)J\right] = 0.
\end{equation}

Therefore, $\omega\left(\frac{\dot{\phi}}{\omega}J + \dot{J}\right)+(\beta - a\cos T)J=0$, so that 
\begin{equation}
    \dot{\phi}\left(J_0 + \omega J_1 + \ldots\right) + \omega\left(\dot{J}_0 + \omega \dot{J}_1 + \ldots\right) + \left(\beta - a\cos T\right)\left(J_0 + \omega J_1 + \ldots\right) = 0.
\end{equation}

The $\mathcal{O}(1)$ corresponding equation is 
\begin{equation}
    J_0\left(\dot{\phi} + (\beta - a \cos T)\right) = 0,
\end{equation}
and solving it for $\phi$ gives 
\begin{equation}
    \phi(T) = \phi^*  -\beta T + a \sin T.
\end{equation}

The $\mathcal{O}(\omega)$ corresponding equation is 
\begin{equation}
    J_1\left(\dot{\phi} + (\beta - a \cos T)\right) + \dot{J}_0 = 0,
\end{equation}
so, as $\dot{\phi} + (\beta - a \cos T)=0$, we obtain $J_0 = J_0^*$ (constant).

Consequently, Eq.~(\ref{eq:WKBJ}) becomes 
\begin{equation}
    y = K\exp\left(\frac{-\beta T+a\sin T}{\omega}\right),
\end{equation}
for $K = p\exp(\frac{\phi^*}{\omega})J_0$ constant. Using that $T = \omega t$ and the initial value $y(0) = y^*$, we obtain the approximation 
\begin{equation}
    y = \left(r^* - \frac{\beta}{\beta - a} \right)\exp\left(-\beta t +\frac{a}{\omega}\sin(\omega t)\right).
\end{equation}

Finally, as $r = y + z(\omega t)$, we have 
\begin{equation}
   r \sim \frac{\beta}{\beta -a\cos(\omega t)} + \left(r^* - \frac{\beta}{\beta - a} \right)\exp\left(-\beta t +\frac{a}{\omega}\sin(\omega t)\right).
\end{equation}

\subsection{Fast variations of the gradients, \texorpdfstring{$\bm{\omega \gg 1}$}.}

We solve now the problem using the method of multiple scales. Let us assume that $r(t) = r(T_1,T_2)$, where $T_1 = t$ and $T_2 = \omega t$, so that 
\begin{equation}
    \frac{dr}{dt} = \frac{\partial r}{\partial T_1}\frac{\partial T_1}{\partial t} + \frac{\partial r}{\partial T_2}\frac{\partial T_2}{\partial t} = \frac{\partial r}{\partial T_1} + \omega\frac{\partial r}{\partial T_2}.
\end{equation}

If $\varepsilon = \frac{1}{\omega} \ll 1$, Eq.~(\ref{eq:a_ini}) becomes 
\begin{equation}
    \varepsilon\left(\frac{\partial r}{\partial T_1} + \frac{\partial r}{\partial T_2}\right) + \varepsilon\left(\beta - a\cos(T_2)\right)r= \varepsilon \beta.
\end{equation}

If we use an asymptotic expansion of $r$, $r = r_0 +\varepsilon r_1$ we obtain:
\begin{equation}
    \varepsilon\left(\frac{\partial r_0}{\partial T_1} + \varepsilon\frac{\partial r_1}{\partial T_1}\right) + \frac{\partial r_0}{\partial T_2} + \varepsilon \frac{\partial r_1}{\partial T_2} +  \varepsilon\left(\beta - a\cos(T_2)\right)\left(r_0 + \varepsilon r_1\right) + \mathcal{O}(\varepsilon^2) = \varepsilon \beta.
\end{equation}

Solving the equation obtained collecting the  $\mathcal{O}(1)$ terms, we find $r_0$ is a function of $T_1$ only, that is 
\begin{equation}
    r_0 = r_0(T_1), 
\end{equation}
with the initial condition $r_0(0) = r^*.$

Now, for the equation obtained collecting the $\mathcal{O}(\varepsilon)$ terms, we have
\begin{equation} \label{eq:condition_1}
    r_0'(T_1) + \frac{\partial r_1}{\partial T_2} + (\beta - a\cos(T_2))r_0(T_1) = \beta.
\end{equation}

As $r_1$ is a periodic correction, integrating Eq.~(\ref{eq:condition_1}) over  $T_2\in [0,2\pi]$ we obtain 
\begin{equation} \label{eq:integral}
    2\pi r_0'(T_1) + 2\pi\beta r_0(T_1) = 2\pi \beta,
\end{equation}
and thus, solving for $r_0(T_1)$, we have 
\begin{equation}
    r_0(T_1) = r^*e^{-\beta T_1} + (1-e^{-\beta T_1}).
\end{equation}

Hence, the leading order approximation is 
\begin{equation}
    r(t) \sim r^*e^{-\beta t} + (1-e^{-\beta t}).
\end{equation}

Now, for the $\mathcal{O}(\varepsilon)$ equation we have:
\begin{equation} \label{eq:oeps}
    \frac{\partial r_1}{\partial T_2} = a\cos(T_2)r_0(T_1).
\end{equation}

Since $r_1(T_1,T_2 = 0) = 0$,
\begin{equation} \label{eq:v1_sol}
    r_1(T_1,T_2) = a\sin(T_2)r_0(T_1).
\end{equation}
and therefore the first order correction is, $r_0 + \varepsilon r_1$, that is:
\begin{equation}
    r(t) \sim \left(r^*e^{-\beta t} + (1-e^{-\beta t})\right)\left(1 + \frac{a}{\omega}\sin(\omega t)\right).
\end{equation}

\subsection{Dominant chemotaxis, \texorpdfstring{$\bm{a\gg \beta.}$}.}\label{sec12}

We solve now the problem using the standard asymptotic expansion method. We set $\varepsilon = \beta/a \ll 1$. Then, Eq.~(\ref{eq:a_ini}) reduces to 
\begin{equation}
    r' + a\left(\varepsilon-\cos(\omega t)\right)r = a\varepsilon.
\end{equation}

The leading order solution is obtained immediately as it is the solution to the homogeneous linear differential equation:
\begin{equation}
    r_0'-a\cos(\omega t)r_0=0.
\end{equation}

Thus  the leading order approximation is 
\begin{equation}
   r\sim  r_0=r^*e^{\frac{1}{\omega}a\sin(\omega t)}, 
\end{equation}
though the next order correction generates a cumbersome expression and thus is not presented.

\subsection{Dominant cell proliferation,  \texorpdfstring{$\bm{a\ll \beta.}$}.}

Again, we solve the problem using the standard asymptotic expansion method. Now, we set $\varepsilon = a/\beta \ll 1$. Then, Eq.~(\ref{eq:a_ini}) gives 
\begin{equation}
    r' + \beta \left(1-\varepsilon\cos(\omega t)\right)r = \beta.
\end{equation}

The leading order solution is obtained immediately as  the solution to the inhomogeneous linear differential equation 
\begin{equation}
    r_0'+\beta r_0=\beta, 
\end{equation}
and hence is given by 
\begin{equation}
    r \sim r_0=r^*e^{-\beta t} + \left(1-e^{-\beta t}\right).
\end{equation}
For the first correction, an asymptotic expansion of $r$ of the form $r=r_0 + \varepsilon r_1$ reveals that  
\begin{equation}
    r_1'+\beta r_1=\beta\cos(\omega t)r_0.
\end{equation}
The solution to this ODE with initial condition $r_1(0) = 0$ is given by 
\begin{equation}
    r_1(t) = \frac{\beta\left[(r^*-1)\gamma^2\sin(\omega t)+\omega^2e^{\beta t}\sin(\omega t) + \omega \beta e^{\beta t}\cos(\omega t)\right]}{\omega \gamma^2e^{\beta t}}-\frac{\beta^2}{\gamma^2}e^{-\beta t},
\end{equation}
where we have defined $\gamma^2 = \beta^2+\omega^2$. Hence up to the first order correction we have 
\begin{equation}
    r(t) \sim r_0 + \varepsilon r_1 =  1+\frac{(r^*-1)\gamma^2 \omega-A\left[\left(\gamma^2r^*-\gamma+\omega^2e^{\beta t}\right)\sin(\omega t)+ \omega \beta e^{\beta t}\cos(\omega t)-\beta\omega\right]}{\omega \gamma^2 e^{\beta t}}.
\end{equation}

\clearpage
 
\section{Simplification of the transport equation for chemotaxis}\label{appft}

With $v(t,x)$ denoting the concentration of the chemoattractant, and with $\Pi_3\gg 1$ we have 
in section \ref{tgm} an equation of the form 
$$ v_t = \Pi_3 v_{xx}, $$ 
since $w(u,v) =0$ has been imposed on \eqref{eq:governing_dimensionless_1b}, while in section \ref{ccoc} we have equations of one of the two forms
$$ v_t = \Pi_3 v_{xx} - \Pi_3\lambda , ~~~~~~~ v_t = \Pi_3 v_{xx}
-\Pi_3 \lambda v, $$ via Equations (\ref{vint1}), (\ref{vint2}) respectively.
These are accompanied by boundary conditions of the form 
$$v(t,0)= \psi_1(t), ~~~~~~~v(t,1) = \psi_2(t), $$ 
and initial conditions $v(0,x)=v_0(x)$ are required to close the system. For simplicity, we assume the boundary conditions and initial conditions are consistent at $(t,x)=(0,0), (0,1)$.

To further proceed we firstly assume  $\Pi_3^{-1}\ll \lambda \ll \Pi_3$, so that $\lambda$ can be treated as unit order of magnitude in asymptotic methods based on $\varepsilon = \Pi_3^{-1}\ll 1$, and we also assume that $\psi_1(t),$ $\psi_2(t)$ have derivatives that are unit order of magnitude, or less. With these weak assumptions, our objective is show that the time derivative $v_t$ can be neglected at leading order, justifying the use of Eqs. \eqref{eq:governing_dimensionless_2b}, \eqref{vfin1}, \eqref{vfin2} in the main text and also justifying the neglect of the consideration of initial conditions in the main text on the grounds this only governs fast initial transients. 

Below we work with the PDE 
\begin{eqnarray} \label{vgen} v_t = \Pi_3 v_{xx} - \mu \Pi_3 \lambda v^\zeta, \end{eqnarray} 
with  $\mu,~\zeta\in\{0,1\}$ so that $\mu=0$ gives one the above equations while $\mu=1,~\zeta=0$, gives another with the final possibility corresponding to $\mu=1,~\zeta=1$, allowing  the three cases to be considered  simultaneously. 

We have an outer timescale of $t$ and an inner timescale of $\tau = \Pi_3 t=t/\varepsilon$. In the outer region, $t\gg\epsilon$, for the leading order outer solution $v^\mathrm{out}(t,x)$ one indeed has 
$$ 0 =  v^\mathrm{out}_{xx} - \mu  \lambda (v^\mathrm{out})^\zeta, ~~~~~ v^\mathrm{out}(t,0)= \psi_1(t), ~~~~~~~v^\mathrm{out}(t,1) = \psi_2(t),$$ 
without consideration of the initial condition. Hence the solutions presented in the main text, for instance 
Eqs. \eqref{vout1}, \eqref{vout2},\eqref{vout3},\eqref{vout4}, are the same as the  leading order outer solutions. With the leading order  inner solution $v^\mathrm{in}(\tau,x)$  we have 
$$  v^\mathrm{in}_\tau =  v^\mathrm{in}_{xx} - \mu  \lambda ( v^\mathrm{in})^\zeta, ~~~~ v^\mathrm{in}(0,x)=v_0(x),$$ 
and the  leading order boundary conditions 
$$ v^\mathrm{in}(\tau,0)= \psi_1(\varepsilon \tau) = \psi_1(0) , ~~~~~~
v^\mathrm{in}(\tau,1)= \psi_2(\varepsilon \tau) = \psi_2(0) ,$$ 
on noting $t=\varepsilon \tau$ and where $\mathcal{O}(\varepsilon)$ corrections are dropped in the final term for both boundary conditions. 
Recalling   $\mu,~\zeta\in\{0,1\}$ and using the decomposition 
$$ v^\mathrm{in}(\tau,x) = v^\mathrm{out}(0,x)  + q(\tau,x), $$
without loss of generality, we have
$$ q_\tau = q_{xx} -\mu\zeta \lambda q, ~~~ q(\tau,0)=q(\tau,1)=0, ~~~ q(0,x) = v_0(x)-v^\mathrm{out}(0,x).$$
With the Fourier decomposition of $q(\tau,x)$ and its initial condition 
$$ q(x,\tau) = \sum_{n=1}^\infty q_n(\tau) \sin\left(n\pi x \right), ~~~ q(0,x) = \sum_{n=1}^\infty q_n^0 \sin\left(n\pi x \right), 
$$
we have $$ \frac{\mathrm{d}q_n}{\mathrm{d}\tau} = -n^2\pi^2 q_n - \mu \zeta\lambda q_n, ~~~ q_n(0) = q_n^0, $$ so that 
$$  q(x,\tau) = \sum_{n=1}^\infty q_n^0 \sin\left(n\pi x \right) \mathrm{e}^{-(n^2\pi^2+\mu\zeta\lambda)\tau} \approx 0 ~~ \mbox{for} ~~ \pi^2 \tau \gg 1.$$
Thus the leading order composite solution is given by 
\begin{eqnarray} v(t,x) &=& v^\mathrm{out}(t,x) + v^\mathrm{in}(\tau(t),x) - \lim_{t \rightarrow 0}v^\mathrm{out}(t,x) =v^\mathrm{out}(t,x) + q(\tau(t),x)
 \nonumber \\  &\approx & v^\mathrm{out}(t,x) ~~ \mbox{for} ~~ t  \gg \epsilon/\pi^2.
\end{eqnarray}
Thus, as implemented in the main text, working solely  with the outer solution and neglecting  the initial conditions is a rational asymptotic approximation at leading order with respect to $\varepsilon = \Pi_3^{-1} \ll 1$, once initial transients have decayed, that is for times satisfying $t\gg \epsilon/\pi^2 = \Pi_3^{-1}/\pi^2$. As the initial transients, and very short times, are not of interest for determining the behaviour of the invasive front of cells for the majority of its propagation we thus work only with the outer equations and solution in the main text.

\end{document}